\RequirePackage{fix-cm}
\documentclass[twocolumn]{svjour3}          
\smartqed  
\usepackage{graphicx}
\pdfoutput=1

\graphicspath{{./}{./figs/}}

  \def\clap#1{\hbox to 0pt{\hss#1\hss}}

\providecommand{\mat}[1]{\bm{#1}}%
\renewcommand{\vec}[1]{\mathbf{#1}}
\newcommand{\vecalt}[1]{\bm{#1}}

\providecommand{\mA}{\ensuremath{\mat{A}}}
\providecommand{\mB}{\ensuremath{\mat{B}}}
\providecommand{\mC}{\ensuremath{\mat{C}}}

\providecommand{\mI}{\ensuremath{\mat{I}}}

\providecommand{\mW}{\ensuremath{\mat{W}}}

\providecommand{\mLambda}{\ensuremath{\mat{\Lambda}}}
\providecommand{\mSigma}{\ensuremath{\mat{\Sigma}}}

\providecommand{\mzero}{\ensuremath{\mat{0}}}

\providecommand{\va}{\ensuremath{\vec{a}}}
\providecommand{\vb}{\ensuremath{\vec{b}}}

\providecommand{\vw}{\ensuremath{\vec{w}}}
\providecommand{\vx}{\ensuremath{\vec{x}}}

\providecommand{\vz}{\ensuremath{\vec{z}}}

\providecommand{\vmu}{\ensuremath{\vecalt{\mu}}}

\newcommand{\hvw}{\hat{\vw}}

\newcommand{\sI}{\mathcal{I}}

\newcommand{\sO}{\mathcal{O}}

\newcommand{\sS}{\mathcal{S}}

\newcommand{\Exp}[1]{\mathbb{E}\left[#1\right]}

\newcommand{\Cov}[1]{\mathbb{C}\operatorname{ov}\left[#1\right]}

\newcommand{\Var}[1]{\operatorname{Var}\left[#1\right]}

\newcommand{\indep}{\rotatebox[origin=c]{90}{$\models$}}
\newcommand{\bmat}[1]{\begin{bmatrix}#1\end{bmatrix}}

\newcommand{\dist}[2]{\mathrm{dist}\,\left(\colspan{#1},\,\colspan{#2}\right)}

\renewcommand{\span}{\operatorname{span}}

\newcommand{\colspan}[1]{\operatorname{colspan}(#1)}
\newcommand{\nullsp}[1]{\operatorname{null}(#1)}

\newcommand{\SDRS}{\mathcal{S}_{\operatorname{DRS}}}

\newcommand{\CIR}{\mC_{\operatorname{IR}}}

\newcommand{\CSIR}{\mC_{\operatorname{SIR}}}
\newcommand{\hCSIR}{\hat{\mC}_{\operatorname{SIR}}}
\newcommand{\CAVE}{\mC_{\operatorname{AVE}}}

\newcommand{\CSAVE}{\mC_{\operatorname{SAVE}}}
\newcommand{\hCSAVE}{\hat{\mC}_{\operatorname{SAVE}}}

\newcommand{\Bind}{B_{\text{ind}}}
\newcommand{\Nrmin}{N_{r_{\min}}}

\usepackage{psfrag, epsf, epstopdf}
\usepackage{enumerate, enumitem}
\usepackage{array}
\usepackage[hyphens]{url}
\usepackage{hyperref}
\usepackage{csquotes}
\usepackage{natbib}
\usepackage{algorithm, algorithmic}
\usepackage{amsmath, amsfonts, amsthm, mathtools}
\usepackage{tabulary, booktabs}
\usepackage{subfig}
\usepackage{bm}
\usepackage{color}

\newtheorem{thm}{Theorem}

\newtheorem{define}{Definition}

\newtheorem{prob}{Problem}
\newtheorem{assumption}{Assumption}
\newtheorem{ex}{Example}

\begin{document}

\title{Inverse regression for ridge recovery
\thanks{Glaws' work is supported by the Ben L. Fryrear Ph.D. Fellowship in Computational Science at the Colorado School of Mines and the Department of Defense, Defense Advanced Research Projects Agency's program Enabling Quantification of Uncertainty in Physical Systems. Constantine's work is supported by the US Department of Energy Office of Science, Office of Advanced Scientific Computing Research, Applied Mathematics program under Award Number DE-SC-0011077.}
}
\subtitle{A data-driven approach for parameter reduction in computer experiments}


\author{Andrew Glaws \and Paul G. Constantine \and R. Dennis Cook}


\institute{A. Glaws \at
			  Department of Computer Science \\
			  University of Colorado Boulder \\
              \email{andrew.glaws@colorado.edu}
           \and
           P. G. Constantine \at
			  Department of Computer Science \\
			  University of Colorado Boulder \\
              \email{paul.constantine@colorado.edu}   
           \and
           R. D. Cook \at
			  Department of Applied Statistics \\
			  University of Minnesota \\
              \email{dennis@stat.umn.edu}              
}

\date{Received: date / Accepted: date}

\maketitle

\begin{abstract}
Parameter reduction can enable otherwise infeasible design and uncertainty studies with modern computational science models that contain several input parameters. In statistical regression, techniques for \emph{sufficient dimension reduction} (SDR) use data to reduce the predictor dimension of a regression problem. A computational scientist hoping to use SDR for parameter reduction encounters a problem: a computer prediction is best represented by a deterministic function of the inputs, so data comprised of computer simulation queries fail to satisfy the SDR assumptions. To address this problem, we interpret SDR methods \emph{sliced inverse regression} (SIR) and \emph{sliced average variance estimation} (SAVE) as estimating the directions of a \emph{ridge function}, which is a composition of a low-dimensional linear transformation with a nonlinear function. Within this interpretation, SIR and SAVE estimate matrices of integrals whose column spaces are contained in the ridge directions' span; we analyze and numerically verify convergence of these column spaces as the number of computer model queries increases. Moreover, we show example functions that are not ridge functions but whose inverse conditional moment matrices are low-rank. Consequently, the computational scientist should beware when using SIR and SAVE for parameter reduction, since SIR and SAVE may mistakenly suggest that truly important directions are unimportant.
\end{abstract}

\section{Introduction and related literature}
\label{sec:intro}

Advances in computing have enabled complex computer simulations of physical systems that complement traditional theory and experiments. We generically model this type of physics-based input/output system as a function $f(\vx)$, where $\vx$ represents a vector of continuous physical inputs and $f$ represents a continuous physical output of interest. It is common to treat $f$ as a deterministic function, since a bug-free computational model produces the same outputs given the same inputs; in other words, there is no random noise in simulation outputs. This perspective is the de facto setup in \emph{computer experiments}~\citep{Sacks89,Koehler96,Santner03}. The setup is similar in the field of \emph{uncertainty quantification}~\citep{SmithUQ2013,SullivanUQ2015,HandbookUQ}; although the inputs may be modeled as random variables to represent aleatory uncertainty, the input/output map is deterministic.

Scientific studies based on deterministic computer simulations may require multiple evaluations of the model at different input values, which may not be feasible when a single model run is computationally expensive. In this case, one may construct a cheap approximation to the model---sometimes referred to as a response surface~\citep{Myers1995,Jones2001}, surrogate model~\citep{Razavi2012,Allaire2010}, metamodel~\citep{Wang2006}, or emulator~\citep{Challenor2012}---that can then be sampled thoroughly. However, such constructions suffer from the \emph{curse of dimensionality}~\citep{Donoho00,Traub1998}; loosely, the number of function queries needed to construct an accurate response surface grows exponentially with the input parameter dimension. In computational science, this curse is exacerbated by the computational cost of each query, which often involves the numerical approximation of a system of partial differential equations. In high dimensions (i.e., several physical input parameters), the number of queries needed to reach error estimate asymptopia (or even construct the response surface) is generally considered infeasible. The computational scientist may attempt to simplify the model by fixing unimportant input parameters at nominal values, thus reducing the dimension of the input space to the point that response surface modeling becomes practical. This concept is generalized by active subspaces~\citep{Constantine15}, which seek out off-axis anisotropy in the function.

In regression modeling~\citep{Weisberg2005}, the given data are predictor/response pairs $\{ [ \, \vx_i^\top \, , \,  y_i \, ] \}$ that are assumed to be independent realizations of a random vector with a joint predictor/response probability density function. In this context, subspace-based dimension reduction goes by the name \emph{sufficient dimension reduction} (SDR)~\citep{Cook98,Li2018,Adragni09}. Techniques for SDR include sliced inverse regression (SIR)~\citep{Li91}, sliced average variance estimation (SAVE)~\citep{Cook91}, ordinary least squares (OLS)~\citep{OLS89}, and principal Hessian directions (pHd)~\citep{Li92}---among several others. These techniques seek a low-dimensional subspace in the predictor space that is sufficient to statistically characterize the relationship between predictors and response. SDR methods are gaining interest for parameter reduction in computational science~\citep{Li16,Zhang2017,Pan2017}. The first use we know applies OLS, SIR, and pHd to a contaminant transport model from Los Alamos National Labs, where the results revealed simple, exploitable relationships between transformed inputs and the computational model's output~\citep[Section 4]{Cook94b}. However, these papers do not carefully address the important mathematical differences that arise when applying a technique developed for statistical regression to study a deterministic function. 

A data set comprised of point queries from a deterministic function (e.g., an experimental design in the physical inputs and associated computer model predictions) differs from the regression data set (independent draws of a joint predictor/response distribution). For example, the former does not admit a joint density on the input/output space. Thus, the data from computational science simulations fail to satisfy the regression assumptions. This failure has important consequences for interpreting \emph{error} and \emph{convergence} in the response surface. We do not explore these consequences in this paper, but the differences between statistical regression and function approximation are essential to our thesis. Instead, we focus on (i) the interpretation of dimension reduction subspaces from SIR or SAVE in the absence of random noise in the functional relationship between inputs and outputs and (ii) the proper meaning of \emph{sufficient dimension reduction} for deterministic functions.

We evaluate and analyze the \emph{inverse regression} methods SIR and SAVE as candidates for gradient-free subspace-based parameter reduction in deterministic functions. SIR- and SAVE-style dimension reduction follows from low-rank structure in the regression's inverse conditional moment matrices. The appropriate way to interpret SIR and SAVE in this context is as methods for \emph{ridge recovery}. A function of several variables is a \emph{ridge function} when it is a composition of a low-dimensional linear transformation with a nonlinear function of the transformed variables~\citep{Pinkus15}---i.e., $f(\vx)=g(\mA^\top\vx)$ where $\mA$ is a tall matrix. The goal of ridge recovery is to estimate the matrix $\mA$ (more precisely, $\mA$'s column space known as the \emph{ridge subspace}) using only point queries of $f$. Recent work by \cite{Fornaiser12} and related work by \cite{Tyagi2014} develop ridge recovery algorithms that exploit a connection between a linear measurement of a gradient vector and the function's directional derivative to estimate $\mA$ with finite differences. \cite{constantine2015sketching} exploit a similar connection to estimate active subspaces with directional finite differences. In contrast, the inverse regression approach for ridge recovery is based on conditional moments of the predictors given the response; we detail how this difference (gradients versus inverse conditional moments) affects the methods' ability to recover the ridge directions. 

\subsection{Main results and discussion}
We develop and interpret sufficient dimension reduction theory in the context of parameter reduction for a deterministic scalar-valued function of several variables. By this interpretation, SDR methods translate naturally as ridge recovery methods; loosely, \emph{the noiseless analog of SDR is ridge recovery}. 
This interpretation illuminates the theoretical limitations of inverse regression for parameter reduction. We show that, if the data-generating function is a ridge function, then the population (i.e., no error from finite sampling) SIR and SAVE subspaces are (possibly strict) subspaces of the function's ridge subspace. Consequently, the rank of the inverse conditional moment matrices is bounded above by the dimension of the ridge subspace. We also show how the inverse conditional moments become integrals, so that the finite-sample SIR and SAVE algorithms can be viewed as numerical integration methods. By choosing input samples independently at random, SIR and SAVE become Monte Carlo methods; we analyze the convergence of SIR and SAVE subspaces as the number $N$ of samples increases, and we derive the expected $\mathcal{O}_p(N^{-1/2})$ rate, where $\mathcal{O}_p$ denotes convergence in probability and the constant depends inversely on associated eigenvalue gap. Moreover, this view enables more efficient numerical integration methods than Monte Carlo for SIR and SAVE~\citep{Glaws2018}.

\subsection{Practical considerations and scope}
Some readers want to know how these results impact practice before (or in lieu of) digesting the mathematical analysis. In practice, the computational scientist does not know whether her model's input/output map is a ridge function, and she would like to use SIR and SAVE as a tool to test for possible dimension reduction; our theory gives some insight into this use case. First, in the idealized case of no finite-sampling (or numerical integration) error, SIR and SAVE can still mislead the practitioner, since there are non-ridge functions whose inverse conditional moment matrices are low-rank. Arguably, these may be pathological cases that do not arise (frequently) in practice, but quantifying this argument requires mathematically characterizing a \emph{typical problem} found in practice. We believe these failure modes should be well known and guarded against. Second, finite-sample estimates (or numerical approximations) of the inverse conditional moment matrices will never have trailing eigenvalues that are exactly zero due to the finite sampling---even when the data-generating function is an exact ridge functions. Therefore, it may be difficult to distinguish between finite-sample noise in the eigenvalues and directions of small-but-nonzero variability in the function. The latter structure indicates anisotropic parameter dependence that is arguably more common in real world functions than exact ridge functions, and it (the anisotropy) can be exploited for ridge approximation (as opposed to ridge recovery). Because of such indistinguishability, we recommend exercising caution when using the estimated SIR and SAVE subspaces to build ridge approximations. Principally, heuristics such as bootstrapping or cross validation might help distinguish between finite sample noise and directions of relatively small variability. Third, we have found that breaking the elliptic symmetry assumption on the input density (see Section \ref{subsec:SIR_for_RR})---e.g., when the input density is uniform on a hyperrectangle---may produce nonzero eigenvalues in the inverse conditional moment matrices even when the data-generating function is an exact ridge function; distinguishing these effects from finite-sample effects or true variability in the data-generating function may be difficult.

\subsection{Paper outline}
The remainder is structured as follows. Section \ref{sec:SDR} reviews (i) the essential theory of sufficient dimension reduction in statistical regression and (ii) the SIR and SAVE algorithms. In Section \ref{sec:ridge_funcs}, we review ridge functions and the ridge recovery problem. Section \ref{sec:SDR_for_RR} translates SDR to deterministic ridge recovery and analyzes the convergence of the SIR and SAVE subspaces. Section \ref{sec:numericalresults} shows numerical examples that support the analyses using (i) two quadratic functions of ten variables and (ii) a simplified model of magnetohydrodynamics with five input parameters. Section \ref{sec:conclusion} contains concluding remarks. To improve readability, all proofs are provided as supplemental material.

\section{Sufficient dimension reduction}
\label{sec:SDR}

We review the essential theory of sufficient dimension reduction (SDR); our notation and development closely follow Cook's \emph{Regression Graphics: Ideas for Studying Regressions through Graphics}~\citep{Cook98}. The theory of SDR provides a framework for subspace-based dimension reduction in statistical regression. A regression problem begins with predictor/response pairs $\{ [ \, \vx_i^\top \, , \,  y_i \, ] \}$, $i=1,\dots,N$, where $\vx_i \in \mathbb{R}^m$ and $y_i \in \mathbb{R}$ denote the vector-valued predictors and scalar-valued response, respectively. These pairs are assumed to be independent realizations from the random vector $[\, \vx^\top \, , \, y \, ]$ with unknown joint probability density function $p_{\vx, y}$. The object of interest in regression is the conditional random variable $y|\vx$; the statistician uses the given predictor/response pairs to estimate statistics of $y|\vx$---e.g., moments or quantiles. SDR searches for a subspace of the predictors that is sufficient to describe $y|\vx$ with statistics derived from the given predictor/response pairs.

The basic tool of SDR is the \emph{dimension reduction subspace} (DRS).
\begin{define}[Dimension reduction subspace \\ \citep{Cook96}] \label{def:DRS}
Consider a random vector $[\, \vx^\top \, , \, y \, ]$, and let $\mA \in \mathbb{R}^{m \times n}$ with $n \leq m$ be such that
\begin{equation} \label{eq:cond_indep_DRS}
y \, \indep \, \vx | \mA^\top \vx,
\end{equation}
where $\indep$ denotes independence of $y$ and $\vx$. A \emph{dimension reduction subspace} $\SDRS$ for $y|\vx$ is
\begin{equation} \label{eq:DRS}
\SDRS \;=\; \SDRS(\mA) \;=\; \colspan{\mA}.
\end{equation}
\end{define}
Equation \eqref{eq:cond_indep_DRS} denotes the conditional independence of the response and predictors given $\mA^\top \vx$. We can define a new random variable $y | \mA^\top \vx$ that is the response conditioned on the $n$-dimensional vector $\mA^\top \vx$. If $\SDRS(\mA)$ is a DRS for $y|\vx$, then the conditional CDF for $y|\vx$ is the conditional CDF for $y|\mA^\top\vx$~\citep[Chapter 6]{Cook98}. If a regression admits a low-dimensional ($n<m$) DRS, then the predictor dimension can be reduced by considering $y | \mA^\top \vx$. Note that the matrix $\mA$ in \eqref{eq:DRS} does not uniquely define a DRS~\citep{Cook94a}. Any matrix with the same column space is a basis for $\SDRS$

A given regression problem may admit multiple DRSs. We require a uniquely-defined DRS to ensure a well-posed SDR problem. To this end, we consider the \emph{central dimension reduction subspace} or simply the \emph{central subspace}. 
\begin{define}[Central subspace \citep{Cook96}] \label{def:central_DRS}
A subspace $\sS_{y|\vx}$ is a \emph{central subspace} for $y|\vx$ if it is a DRS and $\sS_{y|\vx} \subseteq \SDRS$, where $\SDRS$ is any DRS for $y|\vx$.
\end{define}
When the central subspace exists it is the intersection of all other DRSs,
\begin{equation} \label{eq:central_DRS_intersection}
\sS_{y|\vx} \;=\; \bigcap \sS_{\text{DRS}}
\end{equation}
for all $\SDRS$ of $y|\vx$. The intersection in \eqref{eq:central_DRS_intersection} always defines a subspace, but this subspace need not satisfy the conditional independence from Definition \ref{def:DRS}. Therefore, a regression need not admit a central subspace. However, when a central subspace exists, it is the unique DRS of minimum dimension for $y|\vx$~\citep[Chapter 6]{Cook98}. This minimum dimension is referred to as the \emph{structural dimension} of $y|\vx$.

There are a variety of conditions that ensure the existence of $\sS_{y|\vx}$ for a given regression problem. We consider the following condition on the marginal density $p_{\vx}$ of the predictors.
\begin{thm} [\cite{Cook98}] \label{thm:central_DRS_exist}
Suppose that $\sS_1$ and $\sS_2$ are two DRSs for $y|\vx$ where the marginal density of the predictors $p_{\vx} (\vx)$ has support over a convex set. Then, $\sS_1 \cap \sS_2$ is also a DRS. 
\end{thm}
According to Theorem \ref{thm:central_DRS_exist}, if $p_{\vx}$ has support over a convex set, then any intersection of DRSs will be a DRS---i.e., $\sS_{y|\vx}$ from  \eqref{eq:central_DRS_intersection} is a DRS and (hence) the central subspace. In regression practice, this existence condition may be difficult to verify from the given predictor/response pairs. Existence of the central subspace can be proven under other sets of assumptions, including a generalization of Theorem \ref{thm:central_DRS_exist} to M-sets~\citep{Yin08} and another based on location regressions~\citep[Ch. 6]{Cook98}. The existence criteria in Theorem \ref{thm:central_DRS_exist} is the most pertinent when we employ inverse regression for ridge recovery in Section \ref{sec:SDR_for_RR}.

There are two useful properties of the central subspace that enable convenient transformations. The first involves the effect of affine transformations in the predictor space. Let $\vz = \mB \vx + \vb$ for full rank $\mB \in \mathbb{R}^{m \times m}$ and $\vx \in \mathbb{R}^m$. If $\colspan{\mA} = \sS_{y|\vx}$ for some $\mA \in \mathbb{R}^{m \times n}$, then $\colspan{\mB^{-\top} \mA} = \sS_{y|\vz}$~\citep{Cook96}. This allows us to make the following assumption about the predictor space without loss of generality.
\begin{assumption}[Standardized predictors]
\label{assume:standardized_predictors}
Assume that $\vx$ is standardized such that 
\begin{equation}
\Exp{\vx} = \mzero, \qquad \Cov{\vx} = \mI .
\end{equation}
\end{assumption}
The second property involves mappings of the response. Let $h: \mathbb{R} \rightarrow \mathbb{R}$ be a function applied to the responses that produces a new regression problem, $\{ [ \, \vx_i^\top \, , \, h(y_i) \, ] \}$, $i=1,\dots,N$. The central subspace associated with the new regression problem is contained within the original central subspace,
\begin{equation}
\label{eq:cent_sub_mapped_response}
\sS_{h(y)|\vx}
\;\subseteq\;
\sS_{y|\vx} ,
\end{equation}
with equality holding when $h$ is strictly monotonic~\citep{Cook00a}. Equation \eqref{eq:cent_sub_mapped_response} is essential for studying the \emph{slicing}-based algorithms, SIR and SAVE, for estimating the central subspace, where the mapping $h$ partitions the response space; see Sections \ref{subsec:SIR} and \ref{subsec:SAVE}.

The goal of SDR is to estimate the central subspace for the regression from the given response/predictor pairs.
\begin{prob}[SDR problem] \label{prob:SDR}
Given response/predictor pairs $\{ [ \, \vx_i^\top \, , \, y_i \, ] \}$, with $i=1,\dots,N$, assumed to be independent draws from a random vector $[\, \vx^\top \, , \, y \, ]$ with joint density $p_{\vx,y}$, compute a basis $\mA \in \mathbb{R}^{m \times n}$ for the central subspace $\sS_{y|\vx}$ of the random variable $y|\vx$.
\end{prob}
Next, we review two algorithms for Problem \ref{prob:SDR}---sliced inverse regression~\citep{Li91} and sliced average variance estimation~\citep{Cook91}. These algorithms use the given regression data to approximate population moment matrices of the \emph{inverse regression} $\vx|y$.

\subsection{Sliced inverse regression}
\label{subsec:SIR}

Sliced inverse regression (SIR)~\citep{Li91} is an algorithm for approximating the matrix
\begin{equation}
\label{eq:CIR}
\CIR \;=\; \Cov{\Exp{\vx|y}} 
\end{equation}
using the given predictor/response pairs. $\CIR$ is defined by the \emph{inverse regression function} $\Exp{\vx|y}$, which draws a curve through the $m$-dimensional predictor space parameterized by the scalar-valued response. If the given regression problem admits an elliptically symmetric marginal density $p_{\vx}$, then it satisfies the \emph{linearity condition}~\citep{Eaton86}. Under the linearity condition,
\begin{equation}
\label{eq:CIR_in_Syx}
\colspan{\CIR} \subseteq \sS_{y|\vx} ,
\end{equation}
which means the column space of $\CIR$ will uncover at least part of the central subspace for the regression~\citep{Li91}.

SIR approximates $\CIR$ using given predictor/response pairs by partitioning (i.e., \emph{slicing}) the response space. Consider a partition of the observed response space,
\begin{equation}
\label{eq:slices}
\min_{1 \leq i \leq N} y_i = \tilde{y}_0 < \tilde{y}_1 < \dots < \tilde{y}_{R-1} < \tilde{y}_{R} = \max_{1 \leq i \leq N} y_i ,
\end{equation}
and let $J_r = [\tilde{y}_{r-1}, \tilde{y}_r]$ denote the $r$th  partition for $r = 1,\dots,R$. Define the function
\begin{equation}
\label{eq:hy_r}
h(y) \;=\; r
\quad \text{for} \quad
y \in J_r .
\end{equation}
Applying $h$ to the given responses creates a new regression problem $\{ [ \, \vx_i^\top \, , \, h(y_i) \, ] \}$, where \eqref{eq:CIR} becomes
\begin{equation}
\CSIR
\;=\;
\Cov{\Exp{\vx|h(y)}} .
\end{equation}
The notation $\CSIR$ emphasizes that this matrix is with respect to the sliced version of the original regression problem. Combining \eqref{eq:CIR_in_Syx} and \eqref{eq:cent_sub_mapped_response},
\begin{equation}
\label{eq:SIR_sub_contain}
\colspan{\CSIR}
\;\subseteq\;
\sS_{h(y)|\vx} 
\;\subseteq\;
\sS_{y|\vx} .
\end{equation}
The sliced partition of the response space and the sliced mapping $h$ bin the predictor/response pairs to enable sample estimates of $\Exp{\vx|r}$ for $r = 1,\dots,R$. This is the basic idea behind the SIR algorithm; see Algorithm \ref{alg:SIR}. Note that if the response is discrete, then Algorithm \ref{alg:SIR} produces the maximum likelihood estimator of the central subspace~\citep{Cook09}.

\begin{algorithm}
\caption{Sliced inverse regression \citep{Li91}} \label{alg:SIR}
\textbf{Given:} $N$ samples $\{ [ \, \vx_i^\top \, , \, y_i \, ] \}$, $i = 1,\dots, N$, drawn independently according to $p_{\vx,y}$. \\
\begin{enumerate}
\item Partition the response space as in \eqref{eq:slices}, and let $J_r = [\tilde{y}_{r-1}, \tilde{y}_r]$ for $r = 1,\dots,R$. Let $\sI_r \subset \{ 1,\dots,N \}$ be the set of indices $i$ for which $y_i \in J_r$ and define $N_r$ to be the cardinality of $\sI_r$.
\item For $r = 1,\dots,R$, compute the sample mean $\hat{\vmu}_h (r)$ of the predictors whose associated responses are in $J_r$,
\begin{equation} \label{eq:sample_slice_average}
\hat{\vmu}_h (r) = \frac{1}{N_r} \sum_{i \in \sI_r} \vx_i .
\end{equation}
\item Compute the weighted sample covariance matrix
\begin{equation}
\label{eq:hDSIR}
\hCSIR \;=\; \frac{1}{N} \sum_{r=1}^R N_r \, \hat{\vmu}_h (r) \hat{\vmu}_h (r)^\top .
\end{equation}
\item Compute the eigendecomposition,
\begin{equation}
\hCSIR \;=\; \hat{\mW} \hat{\mLambda} \hat{\mW}^\top,
\end{equation}
where the eigenvalues are in descending order $\hat{\lambda}_1 \geq \hat{\lambda}_2 \geq \dots \geq \hat{\lambda}_m \geq 0$ and the eigenvectors are orthonormal.
\item Let $\hat{\mA} \in \mathbb{R}^{m \times n}$ be the first $n$ eigenvectors of $\hCSIR$; return $\hat{\mA}$.
\end{enumerate}
\end{algorithm} 

Eigenvectors of $\CSIR$ associated with nonzero eigenvalues provide a basis for the SIR subspace, $\colspan{\CSIR}$. If the approximated eigenvalues $\hat{\lambda}_{n+1},\dots,\hat{\lambda}_m$ from Algorithm \ref{alg:SIR} are small, then $m \times n$ matrix $\hat{\mA}$ approximates a basis for this subspace. However, determining the appropriate value of $n$ requires care. \cite{Li91} and \cite{Cook91} propose significance tests based on the distribution of the average of the $m-n$ trailing estimated eigenvalues. These testing methods also apply to the SAVE algorithm in Section \ref{subsec:SAVE}.

For a fixed number of slices, SIR has been shown to be $N^{-1/2}$-consistent for approximating $\colspan{\CSIR}$~\citep{Li91}. In principle, increasing the number of slices may provide improved estimation of the central DRS. However, in practice, SIR's success is relatively insensitive to the number of slices. The number of slices should be chosen such that there are enough samples in each slice to estimate the conditional expectations accurately. For this reason, \cite{Li91} suggests constructing slices such that the response samples are distributed nearly equally.

\subsection{Sliced average variance estimation}
\label{subsec:SAVE}

SAVE uses the variance of the inverse regression. \cite{Li91} recognized the potential for $\Cov{\vx|y}$ to provide insights into the central subspace by noting that, under Assumption \ref{assume:standardized_predictors},
\begin{equation}
\begin{aligned}
\Exp{\Cov{\vx|y}} \;&=\; \Cov{\vx} - \Cov{\Exp{\vx|y}} \\
\;&=\; \mI - \Cov{\Exp{\vx|y}} ,
\end{aligned}
\end{equation}
which can be rewritten as
\begin{equation} \label{eq:cov_identity}
\begin{aligned}
\Cov{\Exp{\vx|y}} \;&=\; \mI - \Exp{\Cov{\vx|y}} \\
\;&=\; \Exp{\mI - \Cov{\vx|y}} .
\end{aligned}
\end{equation}
The left side of \eqref{eq:cov_identity} is precisely $\CIR$. This suggests that $\Exp{\mI - \Cov{\vx|y}}$ may be useful in addressing Problem \ref{prob:SDR}. \cite{Cook91} suggest using $\Exp{(\mI - \Cov{\vx|y})^2}$, which has nonnegative eigenvalues. Define
\begin{equation}
\label{eq:CAVE}
\CAVE
\;=\;
\Exp{\left( \mI - \Cov{\vx|y} \right)^2} .
\end{equation}
Under the \emph{linearity condition}~\citep{Li91} and the \emph{constant covariance condition}~\citep{Cook00a}, the column space of $\CAVE$ is contained within the central subspace,
\begin{equation}
\label{eq:CAVE_in_Syx}
\colspan{\CAVE} \subseteq \sS_{y|\vx} .
\end{equation}
Both of these conditions are satisfied when the predictors have an elliptically symmetric marginal density~\citep{Eaton86}.

Using the partition \eqref{eq:slices} and the map \eqref{eq:hy_r},
\begin{equation}
\label{eq:CSAVE}
\CSAVE
\;=\;
\Exp{\left( \mI - \Cov{\vx|h(y)} \right)^2} .
\end{equation}
The notation $\CSAVE$ indicates application to the sliced version of the original regression problem. Combining \eqref{eq:CAVE_in_Syx} and \eqref{eq:cent_sub_mapped_response},
\begin{equation}
\colspan{\CSAVE}
\;\subseteq\;
\sS_{h(y)|\vx}
\;\subseteq\;
\sS_{y|\vx} .
\end{equation}
Algorithm \ref{alg:SAVE} shows the SAVE algorithm, which computes a basis for the column span of the sample estimate $\hCSAVE$. This basis is a $N^{-1/2}$-consistent estimate of $\colspan{\CSAVE}$~\citep{Cook00a}. Increasing the number of slices improves the estimate but suffers from the same drawbacks as SIR. In practice, SAVE performs poorly compared to SIR when few predictor/response pairs are available. This is due to difficulties approximating the covariances within the slices using too few samples. For this reason, \cite{Cook09} suggest trying both methods to approximate the central DRS.

\begin{algorithm}[]
\caption{Sliced average variance estimation \citep{Cook00a}} \label{alg:SAVE}
\textbf{Given:} $N$ samples $\{ [ \, \vx_i^\top \, , \, y_i \, ] \}$, $i = 1,\dots, N$, drawn independently according to $p_{\vx,y}$. \\
\begin{enumerate}
\item Define a partition of the response space as in \eqref{eq:slices}, and let $J_r = [\tilde{y}_{r-1}, \tilde{y}_r]$ for $r = 1,\dots,R$. Let $\sI_r \subset \{ 1,\dots,N \}$ be the set of indices $i$ for which $y_i \in J_r$ and define $N_r$ to be the cardinality of $\sI_r$.
\item For $r = 1,\dots,R$,
\begin{enumerate}
\item Compute the sample mean $\hat{\vmu}_h (r)$ of the predictors whose associated responses are in the $J_r$,
\begin{equation}
\hat{\vmu}_h (r) = \frac{1}{N_r} \sum_{i \in \sI_r} \vx_i .
\end{equation}
\item Compute the sample covariance $\hat{\mSigma}_h (r)$ of the predictors whose associated responses are in $J_r$,
\begin{equation}
\hat{\mSigma}_h (r) = \frac{1}{N_r - 1} \sum_{i \in \sI_r} \left( \vx_i - \hat{\vmu}_h (r) \right) \left( \vx_i - \hat{\vmu}_h (r) \right)^\top
\end{equation}
\end{enumerate}
\item Compute the matrix,
\begin{equation}
\label{eq:hDSAVE}
\hCSAVE \;=\; \frac{1}{N} \sum_{r=1}^R N_r \left( \mI - \hat{\mSigma}_h (r) \right)^2 .
\end{equation}
\item Compute the eigendecomposition,
\begin{equation}
\hCSAVE \;=\; \hat{\mW} \hat{\mLambda} \hat{\mW}^\top,
\end{equation}
where the eigenvalues are in descending order $\hat{\lambda}_1 \geq \hat{\lambda}_2 \geq \dots \geq \hat{\lambda}_m \geq 0$ and the eigenvectors are orthonormal.
\item Let $\hat{\mA} \in \mathbb{R}^{m \times n}$ be the first $n$ eigenvectors of $\hCSAVE$; return $\hat{\mA}$.
\end{enumerate}
\end{algorithm}
\section{Ridge functions}
\label{sec:ridge_funcs}

The input/output map of a computer simulation is best modeled by a deterministic function,
\begin{equation}
\label{eq:y_fx}
y \;=\; f(\vx) ,
\qquad
\vx \in \mathbb{R}^m ,
\quad
y \in \mathbb{R} .
\end{equation}
We assume the domain of $f$ is equipped with a known input probability measure $\pi_{\vx}$. In an uncertainty quantification context, $\pi_{\vx}$ encodes uncertainty in the physical inputs due to, e.g., experimental measurement error. We assume $\pi_{\vx}$ admits a density function $p_{\vx}: \mathbb{R}^m \rightarrow \mathbb{R}^+$ and that the $m$ inputs are independent. Common choices for $p_{\vx}$ include a multivariate Gaussian and a uniform density over a hyperrectangle. The pair $f$ and $\pi_{\vx}$ induces an unknown push-forward probability measure on the output space, which we denote by $\pi_y$.

Translating sufficient dimension reduction to deterministic functions naturally leads to \emph{ridge functions}~\citep{Pinkus15}.
\begin{define}[Ridge function \citep{Pinkus15}] \label{def:ridge_function}
A function $f:\mathbb{R}^m\rightarrow\mathbb{R}$ is a \emph{ridge function} if there exists a matrix $\mA \in \mathbb{R}^{m \times n}$ with $n < m$ and a function $g: \mathbb{R}^n \rightarrow \mathbb{R}$ such that
\begin{equation}
\label{eq:y_fx_gAx}
y \;=\; f(\vx) \;=\; g(\mA^\top \vx) .
\end{equation}
The columns of $\mA$ are the \emph{directions} of the ridge function and $g: \mathbb{R}^n \rightarrow \mathbb{R}$ is the \emph{ridge profile}.
\end{define}
A ridge function is nominally a function of $m$ inputs, but it is intrinsically a function of $n < m$ derived inputs. The term \emph{ridge function} sometimes refers to the $n = 1$ case, and $n>1$ is called a \emph{generalized ridge functions}~\citep{Pinkus15}. We do not distinguish between the $n = 1$ and $n > 1$ cases and refer to all such functions as ridge functions for convenience.

Ridge functions appear in a wide range of computational techniques. For example, projection pursuit regression uses a regression model that is a sum of $n$ one-dimensional ridge functions, $\sum_{i=1}^n g_i (\va_i^\top \vx)$, where each $g_i$ is a spline or other nonparametric model~\citep{Friedman80}. Neural network nodes use functions of the form $\sigma (\mW^\top \vx + \vb)$, where $\mW$ is a matrix of weights from the previous layer of neurons, $\vb$ is a bias term for the model, and $\sigma (\cdot)$ is an activation function~\citep{Goodfellow16}.

Ridge functions are good low-dimensional models for functions that exhibit off-axis anisotropic dependence on the inputs since they are constant along the $m-n$ directions of the input space that are orthogonal to the columns of $\mA$. Consider a vector $\vw \in \nullsp{\mA^\top}$. Then, for any $\vx \in \mathbb{R}^m$,
\begin{equation}
\begin{aligned}
f(\vx + \vw)
\;&=\; g(\mA^\top (\vx + \vw)) \\
\;&=\; g(\mA^\top \vx + \mzero) \\
\;&=\; g(\mA^\top \vx) \\
\;&=\; f(\vx) .
\end{aligned}
\end{equation}
We refer to the problem of finding the directions of a ridge function as \emph{ridge recovery}.
\begin{prob}[Ridge recovery] \label{prob:ridge_recovery}
Given a query-able deterministic function $f: \mathbb{R}^m \rightarrow \mathbb{R}$ that is assumed to be a ridge function and an input probability measure $\pi_{\vx}$, find $\mA \in \mathbb{R}^{m \times n}$ with $n<m$ such that
\begin{equation} \label{eq:ridge_function}
f(\vx) \;=\; g(\mA^\top \vx) ,
\end{equation}
for some $g:\mathbb{R}^n \rightarrow \mathbb{R}$.
\end{prob}
Recent papers in signal processing and approximation theory propose and analyze algorithms for ridge recovery. \cite{Cohen12} consider the case of $n=1$ with conditions on the ridge direction. \cite{Fornaiser12} propose an algorithm for general $n<m$ under conditions on $\mA$. \cite{Tyagi2014} extend the approximation results from Fornasier et al.~to more general $\mA$'s. 

Note that the ridge recovery problem is distinct from ridge approximation~\citep{constantine2016near, Hokanson17}, where the goal is to find $\mA$ and construct $g$ that minimize the approximation error for a given $f$.\footnote{Regarding nomenclature, neither ridge recovery nor ridge approximation is related to ridge regression, where Tikhonov regularization is applied to the regression model coefficients~\citep[Chapter 3.4]{Hastie2009}.} Inverse regression may be useful for ridge approximation or identifying near-ridge structure in a given function, but pursuing these ideas is beyond the scope of this manuscript. \cite{Li16} distinguish between an ``ideal scenario''---which loosely corresponds to ridge recovery and where they can prove certain convergence results for inverse regression---and the ``general setting''---where they apply inverse regression as a heuristic in numerical examples.

\section{Inverse regression as ridge recovery}
\label{sec:SDR_for_RR}

This section develops  SIR (Algoroithm \ref{alg:SIR}) and SAVE (Algorithm \ref{alg:SAVE}) as tools for ridge recovery (Problem \ref{prob:ridge_recovery}). The next theorem connects the dimension reduction subspace (Definition \ref{def:DRS}) to the ridge directions (Definition \ref{def:ridge_function}).
\begin{thm}
\label{thm:CondIndep_RidgeFunc}
Let $(\Omega, \Sigma, P)$ be a probability triple. 
Suppose that $\vx: \Omega \rightarrow \mathbb{R}^m$ and $y: \Omega \rightarrow \mathbb{R}$ are random variables related by a measurable function $f: \mathbb{R}^m \rightarrow \mathbb{R}$ so that $y = f(\vx)$. Let $\mA \in \mathbb{R}^{m \times n}$ be a constant matrix. Then $y \, \indep \, \vx | \mA^\top \vx$ if and only if $y = g(\mA^\top \vx)$ where $g: \mathbb{R}^n \rightarrow \mathbb{R}$ is a measurable function.
\end{thm}
Theorem \ref{thm:CondIndep_RidgeFunc} states that conditional independence of the inputs and output in a deterministic function is equivalent to the function being a ridge function. This provides a subspace-based perspective on ridge functions that uses the DRS as the foundation. That is, the directions of the ridge function are relatively unimportant compared to the subspace they span when capturing ridge structure of $f$ with sufficient dimension reduction.

The central subspace (Definition \ref{def:central_DRS})  corresponds to the unique subspace of smallest dimension that completely describes the ridge structure of $f$. Our focus on the subspace instead of the precise basis implies that we can assume standardized inputs $\vx$ (as in Assumption \ref{assume:standardized_predictors}) without loss of generality, which simplifies discussion of the inverse conditional moment matrices that define the SIR and SAVE subspaces for ridge recovery.

Theorem \ref{thm:central_DRS_exist} guarantees existence of the central subspace in regression when the marginal density of the predictors has convex support. However, this condition is typically difficult or impossible to verify for the regression problem. In contrast, the deterministic function is accompanied by a known input measure $\pi_{\vx}$ with density function $p_\vx$. In practice, a relatively non-informative measure is used such as a multivariate Gaussian or uniform density on a hyper-rectangle defined by the ranges of each physical input parameter. Such choices for modeling input uncertainty satisfy Theorem \ref{thm:central_DRS_exist} and guarantee the existence of the central subspace.

The input probability measure can influence the structure of the central subspace. For example, let $\va, \vb \in \mathbb{R}^m$ be constant vectors pointing in different directions and consider the function
\begin{equation}
y
\;=\;
f(\vx)
\;=\;
\left\{ \begin{array}{ll}
(\va^\top \vx)^2 & \text{if } \vx_i > 0 \text{ for } i = 1,\dots,m \\
(\vb^\top \vx)^2 & \text{otherwise.}
\end{array} \right.
\end{equation}
If $\pi_{\vx}$ has a density function with support only for positive values of $\vx$ (e.g., uniform over $[0,1]^m$), then the central subspace is $\span\{\va\}$. Alternatively, if the input density has support over all of $\mathbb{R}^m$ (e.g., multivariate Gaussian), then the central subspace is $\span \{\va, \vb\}$. For this reason, we write the central subspace for a deterministic function as $\sS_{f, \pi_{\vx}}$ to emphasize that this subspace is a property of the given function and input probability measure.

For deterministic functions, the translation of the inverse regression $\vx|y$ (see Section \ref{sec:SDR}) is the inverse image of $f$ for the output value $y$,
\begin{equation}
f^{-1} (y)
\;=\;
\left\{ \, \vx \in \mathbb{R}^m \, : \, f(\vx) = y \, \right\} .
\end{equation}
Unlike the inverse regression $\vx|y$ (which is a random vector), $f^{-1}$ is a fixed set determined by $f$'s contours. Furthermore, the inverse image begets the conditional probability measure $\pi_{\vx|y}$, which is the restriction of $\pi_{\vx}$ to the set $f^{-1}(y)$~\citep{Chang1997}.

\subsection{Sliced inverse regression for ridge recovery}
\label{subsec:SIR_for_RR}

For deterministic functions, we can write $\CIR$ from \eqref{eq:CIR} as an integral:
\begin{equation}
\label{eq:CIR_det}
\CIR
\;=\;
\int \vmu(y) \, \vmu(y)^\top \, d \pi_y (y) 
\end{equation}
where the conditional expectation over the inverse image $f^{-1} (y)$ is
\begin{equation}
\label{eq:cond_exp_det}
\vmu(y)
\;=\;
\int \vx \, d \pi_{\vx|y} (\vx) .
\end{equation}
This term represents the average of all input values that map to a fixed value of the output.

To understand how $\CIR$ can be used for dimension reduction in deterministic functions, consider the following. For $\vw \in \mathbb{R}^m$ with unit norm,
\begin{equation} \label{eq:wT_CIR_w}
\begin{aligned}
\vw^\top \CIR \vw
\;&=\; \vw^\top \left( \int  \vmu (y) \, \vmu (y)^\top  \, d\pi_y(y) \right) \vw \\
\;&=\; \int \left( \vmu (y)^\top \vw \right)^2 \, d\pi_y (y) .
\end{aligned}
\end{equation}
If $\vw \in \nullsp{\CIR}$, then one possibility is that $\vmu(y)$ is orthogonal to $\vw$ for all $y$. The following theorem relates this case to possible ridge structure in $f$.
\begin{thm}
\label{thm:det_SIR}
Let $f: \mathbb{R}^m \rightarrow \mathbb{R}$ with input probability measure $\pi_{\vx}$ admit a central subspace $\sS_{f, \pi_{\vx}}$, and assume $\pi_{\vx}$ admits an elliptically symmetric and standardized density function. Then, $\colspan{\CIR} \subseteq \sS_{f, \pi_{\vx}}$.
\end{thm}
Theorem \ref{thm:det_SIR} shows that a basis for the range of $\CIR$ can be used to estimate $f$'s central subspace. However, this idea has two important limitations. First, by \eqref{eq:wT_CIR_w}, we can write the inner product in the rightmost integrand in terms of the cosine of the angle between the vectors $\vmu (y)$ and $\vw$,
\begin{equation}
\label{eq:wT_CIR_w2}
\vw^\top \CIR \vw
\;=\;
\int \left| \left| \vmu (y) \right| \right|_2^2 \cos^2 (\theta (y)) \, d \pi_y (y) ,
\end{equation}
where $\theta (y)$ is the angle between $\vmu(y)$ and $\vw$. Theorem \ref{thm:det_SIR} uses orthogonality of $\vmu (y)$ and $\vw$ (i.e., $\cos (\theta (y)) = 0$ for all $y$) to show containment of the column space of $\CIR$ within the central subspace; however, the integrand in \eqref{eq:wT_CIR_w2} also contains the squared 2-norm of $\vmu(y)$, which does not depend on $\vw$. If $\vmu(y) = \mzero$ for all $y$, then $\vw^\top \CIR \vw = 0$ for all $\vw$. Consider the following example.
\begin{ex}\label{ex:strict_subset}
Assume $\vx \in \mathbb{R}^2$ is weighted with a bivariate standard Gaussian. Let $y = f (\vx) = x_1 x_2$. For any value of $y$, $\vx \in f^{-1} (y)$ implies $-\vx \in f^{-1} (y)$. Therefore, $\vmu (y) = \mzero$ for all $y$, and $\CIR = \mzero$. But $y$ is not constant over $\mathbb{R}^2$. Thus, $\{\mzero\} = \colspan{\CIR} \subset \sS_{f, \pi_{\vx}} = \mathbb{R}^2$.
\end{ex}
This example shows how $\CIR$, as a tool for ridge recovery, can mislead the practitioner by suggesting ridge structure in a function that is not a ridge function. This could lead one to ignore input space directions that should not be ignored. Note that if we shift the function such that $y = f (\vx) = (x_1 + c_1) (x_2 + c_2)$ for some constants $c_1, c_2 \neq 0$, then the symmetry is broken and $\vmu(y) \neq 0$ for all $y$. In this case, $\CIR$ will recover the central subspace (i.e., all of $\mathbb{R}^2$).

The second limitation of $\CIR$ for ridge recovery follows from the required elliptic symmetry of the input density $p_{\vx}$. This assumption is satisfied if $p_{\vx}$ is a multivariate Gaussian, but it is violated when the density is uniform over the $m$-dimension hypercube. If $f$ is a ridge function and $\vw \in \nullsp{\mA^\top}$, then $\vx \in f^{-1} (y)$ implies $\vx + \vw \in f^{-1} (y)$ so that $f^{-1} (y)$ can be expressed as the union of lines parallel to $\vw$. If the inputs are weighted by an elliptically symmetric density, then the expectation over $f^{-1} (y)$ will be centered such that $\vmu(y)$ is orthogonal to $\vw$. If the inputs do not have an elliptically symmetric density, then the weighting can cause the conditional expectation to deviate in the direction of $\vw$. The magnitude of this deviation also depends on the magnitude of the conditional expectation $||\vmu(y)||_2$.

Next, we examine the sliced approximation of $\CIR$. Recall the output partition from \eqref{eq:slices} and the slicing map $h(y)$ \eqref{eq:hy_r}. Applying $h$ to the deterministic function $y = f(\vx)$ produces the discretized function $r = h(y) = h(f(\vx))$, where $r \in \{ 1,\dots,R \}$. The output space of $h$ is weighted by the probability mass function
\begin{equation}
\label{eq:slice_density}
\omega (r) \;=\; \int_{J_r} \, d \pi_y (y) , \qquad r\in\{1,\dots, R\}.
\end{equation}
Without loss of generality, we assume that the slices are constructed such that $\omega (r) > 0$ for all $r$. If $\omega (r) = 0$ for some $r$, then we can combine this slice with an adjacent slice. The conditional expectation for the sliced output is
\begin{equation}
\label{eq:cond_exp_r}
\vmu_h (r)
\;=\;
\int \vx \, d \pi_{\vx|r} (\vx) ,
\end{equation}
where $\pi_{\vx|r}$ is the conditional measure defined over the set $f^{-1} ( h^{-1} (r) ) =  \{ \, \vx \in \mathbb{R}^m \, : \, h(f(\vx)) = h(y) = r \, \}$. Using \eqref{eq:slice_density} and \eqref{eq:cond_exp_r}, the sliced version of $\CIR$ is
\begin{equation}
\label{eq:CSIR}
\CSIR \;=\; \sum_{r=1}^R \omega (r) \, \vmu_h (r) \, \vmu_h (r)^\top .
\end{equation}
By Theorem \ref{thm:CondIndep_RidgeFunc}, properties of the central subspace extend to the ridge recovery problem. This includes containment of the central subspace under any mapping of the output,
\begin{equation}
\label{eq:SIR_sub_contain_det}
\colspan{\CSIR} \;\subseteq\; \sS_{h \circ f, \pi_{\vx}} \;\subseteq\; \sS_{f, \pi_{\vx}} .
\end{equation}
By approximating $\CSIR$, we obtain an approximation of part of $f$'s central subspace. An important corollary of \eqref{eq:SIR_sub_contain_det} is that the rank of $\CSIR$ is bounded above by the dimension of $\sS_{f,\pi_{\vx}}$.

Note that $\CSIR$ from \eqref{eq:CSIR} is a finite sum approximation of the integral in $\CIR$ from \eqref{eq:CIR_det}. Since $f^{-1} ( h^{-1} (r)) = \cup_{y \in J_r} f^{-1} (y)$, then  $\vmu_h (r)$ is the average of the conditional expectations with $y \in J_r$. That is,
\begin{equation}
\vmu_h (r)
\;=\;
\int_{J_r} \vmu (y) \, d \pi_y (y) .
\end{equation}
Therefore, $\CSIR$ approximates $\CIR$ by a weighted sum of the average values of $\vmu(y)$ within each slice. If $\vmu(y)$ is continuous almost everywhere with respect to $\pi_y$, then $\CIR$ is Riemann-integrable~\citep[Ch. 2]{Folland99}. This ensures that sum approximations using the supremum and infimums of $\vmu(y)$ over each slice converge to the same value. By the sandwich theorem, the average value will converge to this value as well~\citep{Abbott01}. Therefore, $\CSIR$ is a Riemann sum approximation of $\CIR$; as the number of slices $R$ increases, $\CSIR$ converges to $\CIR$.

We turn attention to asymptotic convergence of Algorithm \ref{alg:SIR} for ridge recovery. To generate the data for Algorithm \ref{alg:SIR}, we choose $N$ points $\{\vx_i\}$ in the input space consistent with $\pi_{\vx}$. For each $\vx_i$, we query the function to produce the corresponding output $y_i=f(\vx_i)$. In the computational science context, this corresponds to running the simulation model at particular sets of inputs. If we choose each $\vx_i$ independently according to $\pi_{\vx}$, then we can analyze SIR as a Monte Carlo method for estimating $\CSIR$ from \eqref{eq:CSIR}. Given the input/output pairs, Algorithm \ref{alg:SIR} constructs the random matrix $\hCSIR$. To be clear, $\hCSIR$ is a random estimate of $\CSIR$ because of how we chose the points $\{\vx_i\}$---not because of any randomness in the map $f$. Eigenpairs derived from $\hCSIR$ are also random, and the convergence analysis for Algorithm \ref{alg:SIR} is probabilistic.

The convergence depends on the smallest number of samples per slice over all the slices:
\begin{equation} \label{eq:N_r_min}
N_{r_{\min}} \;=\; \min_{1 \leq r \leq R} N_r,
\end{equation}
where $N_r$ is from Algorithm \ref{alg:SIR}. Recall that the slices are assumed to be constructed such that $\omega(r) > 0$. Thus, $N_{r_{\min}} > 0$ with probability 1 as $N \rightarrow \infty$. The following theorem shows that the eigenvalues of $\hCSIR$ converge to those of $\CSIR$ in a mean-squared sense. 
\begin{thm}
\label{thm:SIR_eig_converge}
Assume that Algorithm \ref{alg:SIR} has been applied to the data set $\{ [ \, \vx_i^\top \, , \, y_i \, ] \}$, with $i=1,\dots,N$, where the $\vx_i$ are drawn independently according to $\pi_{\vx}$ and $y_i = f(\vx_i)$ are point evaluations of $f$. Then, for $k = 1,\dots,m$,
\begin{equation}
\Exp{\left( \lambda_k (\CSIR) - \lambda_k (\hCSIR) \right)^2} \;=\; \sO (N^{-1}_{r_{\min}})
\end{equation}
where $\lambda_k (\cdot)$ denotes the $k$th eigenvalue of the given matrix.
\end{thm}
In words, the mean-squared error in the eigenvalues of $\hCSIR$ decays at a $N_{r_{\min}}^{-1}$ rate. Since $\omega (r) > 0$ for all $r$, $N_{r_{\min}} \rightarrow \infty$ as $N \rightarrow \infty$. Moreover, the convergence rate suggests that one should choose the slices in Algorithm \ref{alg:SIR} such that the same number of samples appears in each slice. This maximizes $N_{r_{\min}}$ and reduces the error in the eigenvalues. 

An important consequence of Theorem \ref{thm:SIR_eig_converge} is that the column space of the finite-sample $\hCSIR$ is not contained in $f$'s central subspace. In fact, due to finite sampling, $\hCSIR$ is not precisely low-rank. With a fixed number of samples, it is difficult to distinguish the effects of finite sampling from actual variability in $f$. However, the eigenvalue convergence implies that one can devise practical convergence tests for low-rank-ness based on sets of samples with increasing size.

The next theorem shows the value of understanding the approximation errors in the eigenvalues for quantifying the approximation errors in the subspaces. We measure convergence of the subspace estimates using the subspace distance~\citep[Chapter 2.5]{Golub96},
\begin{equation} \label{eq:sub_dist}
\dist{\mA}{\hat{\mA}} \;=\; 
\left\| \mA \mA^\top - \hat{\mA} \hat{\mA}^\top \right\|_2,
\end{equation}
where $\mA, \hat{\mA}$ are the first $n$ eigenvectors of $\CSIR$ and $\hCSIR$, respectively. The distance metric \eqref{eq:sub_dist} is the principal angle between the subspaces $\colspan{\mA}$ and $\colspan{\hat{\mA}}$.
\begin{thm}
\label{thm:SIR_sub_converge}
Assume the same conditions from Theorem \ref{thm:SIR_eig_converge}, and let $\Delta_n = \lambda_n (\CSIR) - \lambda_{n+1} (\CSIR)$ denote the gap between the $n$th and $(n+1)$th eigenvalues of $\CSIR$. Then, for sufficiently large $N$,
\begin{equation}
\dist{\mA}{\hat{\mA}} \;=\; 
\Delta_n^{-1} \; \sO_p (N^{-1/2}_{r_{\min}}) ,
\end{equation}
where $\mA, \hat{\mA}$ are the first $n$ eigenvectors of $\CSIR$ and $\hCSIR$, respectively.
\end{thm}
The subspace error decays with asymptotic rate $N_{r_{\min}}^{-1/2}$. The more interesting result from Theorem \ref{thm:SIR_sub_converge} is the inverse relationship between the subspace error and the magnitude of the gap between the $n$th and $(n+1)$th eigenvalues. That is, a large gap between eigenvalues suggests a better estimate of the subspace for a fixed number of samples. We do not hide this factor in the $\sO$ notation to emphasize the importance of Theorem \ref{thm:SIR_eig_converge}, which provides insights into the accuracy of the estimated eigenvalues of $\CSIR$. 

\subsection{SAVE for ridge recovery}
\label{subsec:SAVE_for_RR}

Similar to \eqref{eq:CIR_det}, we express $\CAVE$ from \eqref{eq:CAVE} as an integral, 
\begin{equation}
\label{eq:CAVE_det}
\CAVE
\;=\;
\int \left( \mI - \mSigma (y) \right)^2 \, d \pi_y (y) .
\end{equation}
The conditional covariance $\mSigma(y)$ in \eqref{eq:CAVE_det} is an integral against the conditional probability measure $\pi_{\vx|y}$, 
\begin{equation}
\label{eq:cond_cov_det}
\mSigma (y)
\;=\;
\int \left( \vx - \vmu(y) \right) \, \left( \vx - \vmu(y) \right)^\top \, d \pi_{\vx|y} (\vx) .
\end{equation}
To see the relationship between the $\CAVE$ matrix and ridge functions, let $\vw \in \mathbb{R}^m$ with unit norm:
\begin{equation} \label{eq:wT_CAVE_w}
\begin{aligned}
\vw^\top \CAVE \vw 
\;&=\; \vw^\top \left( \int \left( \mI - \mSigma (y) \right)^2 \,d\pi_y(y) \right) \vw \\
\;&=\; \int \left| \left| \left( \mI -\mSigma(y) \right) \vw \right| \right|_2^2 \, d \pi_y(y) .
\end{aligned}
\end{equation}
Equation \eqref{eq:wT_CAVE_w} relates the column space of $\CAVE$ to the ridge structure in $f$. 
\begin{thm}
\label{thm:det_SAVE}
Let $f: \mathbb{R}^m \rightarrow \mathbb{R}$ with input probability measure $\pi_{\vx}$ admit a central subspace $\sS_{f, \pi_{\vx}}$, and assume $\pi_{\vx}$ admits an elliptically symmetric and standardized density function. Then, $\colspan{\CAVE} \subseteq \sS_{f, \pi_{\vx}}$.
\end{thm}
This result shows the usefulness of $\CAVE$ for revealing ridge structure in deterministic functions: by estimating the column space of $\CAVE$, we obtain an estimate of a subspace of $f$'s central subspace. However, $\CAVE$ suffers two similar pitfalls as $\CIR$. First, $\CAVE$ can mislead the practitioner by suggesting ridge structure that does not exist---i.e., when $\colspan{\CAVE} \subset \sS_{f, \pi_{\vx}}$---as the following example illustrates. 
\begin{ex}
\label{ex:strict_subset2}
Assume $\vx \in \mathbb{R}^2$ is weighted by a bivariate standard Gaussian. Let
\begin{equation}
\label{eq:ex2_func}
y \;=\; f(\vx) \;=\;
\begin{cases} 
y_1 & \text{if } ||\vx||_2 \leq r_1 \text{ or } ||\vx||_2 \geq r_2 , \\
y_2 & \text{if } r_1 < ||\vx||_2 < r_2 ,
\end{cases}
\end{equation}
for some $0 < r_1 < r_2$. This functions looks like a bullseye with the central circle and outer ring mapping to $y_1$ and the middle ring mapping to $y_2$. If we choose $r_1$ and $r_2$ appropriately, then we can obtain
\begin{equation}
\label{eq:ident_sigmas}
\mSigma(y_1) \;=\; \mSigma(y_2) \;=\; \mI .
\end{equation}
Note that $\vmu(y_1) = \vmu(y_2) = 0$ for all choices of $r_1$ and $r_2$. Then,
\begin{equation}
\begin{aligned}
\mSigma(y_1) \;&=\; \int \vx \vx^\top \, d \pi_{\vx|y_1} (\vx) , \\
\;&=\; \left( 1 + \frac{r_2^2 e^{-r_2^2/2} - r_1^2 e^{-r_1^2/2}}{2 \left( 1 + e^{-r_2^2/2} - e^{-r_1^2/2} \right)} \right) \bmat{1 & 0 \\ 0 & 1} ,
\end{aligned}
\end{equation}
and
\begin{equation}
\begin{aligned}
\mSigma(y_2) \;&=\; \int \vx \vx^\top \, d \pi_{\vx|y_2} (\vx) , \\
\;&=\; \left( 1 + \frac{r_2^2 e^{-r_2^2/2} - r_1^2 e^{-r_1^2/2}}{2 \left(e^{-r_2^2/2} - e^{-r_1^2/2} \right)} \right) \bmat{1 & 0 \\ 0 & 1} .
\end{aligned}
\end{equation}
Thus, \eqref{eq:ident_sigmas} holds when $r_1^2 e^{-r_1^2/2} = r_2^2 e^{-r_2^2/2}$, provided that $0 < r_1 < r_2$. Choosing $r_1$ and $r_2$ that satisfy these requirements results in 
\begin{equation}
\begin{aligned}
\CAVE \;&=\; \Exp{\left(\mI - \Cov{\vx|y}\right)^2} \\
\;&=\; \Exp{\left(\mI - \mI\right)^2} \\
\;&=\; \Exp{\mzero} \\
\;&=\; \mzero .
\end{aligned}
\end{equation}
However, we can see from inspection of \eqref{eq:ex2_func} that $\sS_{f,\pi_{\vx}} = \mathbb{R}^2$. Thus, $\colspan{\CAVE} \subset \sS_{f,\pi_{\vx}}$.
\end{ex}
Example \ref{ex:strict_subset2} shows one way that $\CAVE$ can be fooled into falsely suggesting low-dimensional structure in a function. Note that the rotational symmetry in $f(\vx)$ is the key feature of this function that results in a degenerate $\CAVE$ matrix. This symmetry also tricks the $\CIR$ matrix. In fact, it can be shown that, in general,
\begin{equation}
\label{eq:CIR_subset_CAVE}
\colspan{\CIR} \subseteq \colspan{\CAVE},
\end{equation}
which suggests that this sort of false positive is less likely to occur with $\CAVE$ than with $\CIR$~\citep[Chap. 5]{Li2018}. The exhaustiveness of $\CAVE$ in capturing central subspace can be proven provided that at least one of $\Exp{\vw^\top \vx | y}$ and $\Var{\vw^\top \vx | y}$ are nondegenerate---i.e., explicitly depends on $y$---for all $\vw \in \sS_{f,\pi_{\vx}}$. Notice that the function described in Example \ref{ex:strict_subset2} violates this exhaustiveness condition.

The second limitation of $\CAVE$ arises from the elliptic symmetry requirement on the input density $p_{\vx}$. When $p_{\vx}$ is not elliptically symmetric, we cannot guarantee that the column space of $\CAVE$ is contained within the central subspace. Thus, a basis for the column space of $\CAVE$ could be contaminated by effects of $\pi_{\vx}$---independent of whether or not $f$ is a ridge function. 

Next, we consider the sliced version of $\CAVE$. We use the same slicing function $h$ from \eqref{eq:hy_r} to approximate $\CAVE$ from \eqref{eq:cond_cov_det}. The sliced approximation of $\CAVE$ is
\begin{equation}
\CSAVE \;=\; \sum_{r=1}^R \omega (r) \, \left( \mI - \mSigma_h (r) \right)^2 ,
\end{equation}
where $\omega (r)$ is the probability mass function from \eqref{eq:slice_density} and 
\begin{equation}
\mSigma_h (r)
\;=\;
\int \left( \vx - \vmu_h(r) \right) \, \left( \vx - \vmu_h(r) \right)^\top \, d \pi_{\vx|r} (\vx) .
\end{equation}
By containment of the central subspace,  
\begin{equation}
\label{eq:CSAVE_contain}
\colspan{\CSAVE} \;\subseteq\; \sS_{h \circ f, \pi_{\vx}} \;\subseteq\; \sS_{f, \pi_{\vx}} .
\end{equation}
We can interpret $\CSAVE$ as a Riemann sum approximation of $\CAVE$ using a similar argument as in Section \ref{subsec:SIR_for_RR}. An important corollary of \eqref{eq:CSAVE_contain} is that the rank of $\CSAVE$ is bounded above by the dimension of $f$'s central subspace. 

Algorithm \ref{alg:SAVE} computes a sample approximation of $\CSAVE$, denoted $\hCSAVE$, using given input/output pairs $\{ [ \, \vx_i^\top \, , \, y_i \, ] \}$. When the $\vx_i$ are sampled independently according to $\pi_{\vx}$ and each $y_i=f(\vx_i)$ is a deterministic function query, we can interpret $\hCSAVE$ as a Monte Carlo approximation to $\CSAVE$. Thus, $\hCSAVE$ and its eigenpairs are random---not because of any randomness in the map $f$ but because of the random choices of $\vx_i$. The following theorem shows the rate of mean-squared convergence of the eigenvalues of $\hCSAVE$. 
\begin{thm}
\label{thm:SAVE_eig_converge}
Assume that Algorithm \ref{alg:SAVE} has been applied to the data set $\{ [ \, \vx_i^\top \, , \, y_i \, ] \}$, with $i=1,\dots,N$, where the $\vx_i$ are drawn independently according to $\pi_{\vx}$ and $y_i = f(\vx_i)$ are point evaluations of $f$. Then, for $k = 1,\dots,m$,
\begin{equation}
\Exp{\left( \lambda_k (\CSAVE) - \lambda_k (\hCSAVE) \right)^2} \;=\; \sO (N^{-1}_{r_{\min}})
\end{equation}
where $\lambda_k (\cdot)$ denotes the $k$th eigenvalue of the given matrix.
\end{thm}
We note that the column space of $\hCSAVE$ is not contained in $f$'s central subspace because of finite sampling. However, using a sequence of estimates with increasing $N$, one may be able to distinguish effects of finite sampling from true directions of variability in $f$.

Next, we examine the convergence of the subspaces Algorithm \ref{alg:SAVE} produces, where the subspace distance is from \eqref{eq:sub_dist}. 
\begin{thm}
\label{thm:SAVE_sub_converge}
Assume the same conditions from Theorem \ref{thm:SAVE_eig_converge}, and let $\Delta_n = \lambda_n (\CSAVE) - \lambda_{n+1} (\CSAVE)$ denote the gap between the $n$th and $(n+1)$th eigenvalues of $\CSAVE$. Then, for sufficiently large $N$,
\begin{equation}
\dist{\mA}{\hat{\mA}} \;=\; 
\Delta_n^{-1} \; \sO_p (N^{-1/2}_{r_{\min}}) ,
\end{equation}
where $\mA, \hat{\mA}$ are the first $n$ eigenvectors of $\CSAVE$ and $\hCSAVE$, respectively.
\end{thm}
The subspace error for Algorithm \ref{alg:SAVE} decays asymptotically like $N_{r_{\min}}^{-1/2}$ with high probability. Similar to the estimated SIR subspace from Algorithm \ref{alg:SIR}, the error depends inversely on the eigenvalue gap. If the gap between the $n$th and $(n+1)$th eigenvalues is large, then the error in the estimated $n$-dimensional subspace is relatively small for a fixed number of samples.

\section{Numerical results}
\label{sec:numericalresults}

We apply SIR and SAVE to three ridge functions to study the methods' applicability for ridge recovery and verify our convergence analysis. Since we are not concerned with using SIR and SAVE for ridge approximation, we do not study the methods' behavior for functions that are not ridge functions. For each ridge function, we provide several graphics including convergence plots of estimated eigenvalues, eigenvalue errors, and subspace errors. One graphic is especially useful for visualizing the structure of the function relative to its central subspace: the sufficient summary plot~\citep{Cook98}. Sufficient summary plots show $y$ versus $\mA^\top \vx$, where $\mA$ has only one or two columns that comprise a basis for the central subspace. Algorithms \ref{alg:SIR} and \ref{alg:SAVE} produce eigenvectors that span an approximation of the SIR and SAVE subspaces, respectively. We use these eigenvectors to construct the low-dimensional inputs for sufficient summary plots. The label \emph{sufficient} is tied to the precise definition of statistical sufficiency in sufficient dimension reduction for regression. 
 
The convergence analysis from Sections \ref{subsec:SIR_for_RR} and \ref{subsec:SAVE_for_RR} assume a fixed slicing of the observed $y$ range. However, Algorithms \ref{alg:SIR} and \ref{alg:SAVE} are implemented using an adaptive slicing approach that attempts to maximize $\Nrmin$ for a given set of data. This is done as a heuristic technique for reducing eigenvalue and subspace errors.

The Python code used to generate the figures is available at \url{https://bitbucket.org/aglaws/inverse-regression-for-ridge-recovery}. The scripts require the dev branch of the Python Active-subspaces Utility Library~\citep{constantine2016python}.

\subsection{One-dimensional quadratic ridge function}
\label{ssec:numex1}
\begin{figure*}[!ht]
\centering
\subfloat[Eigenvalues of $\hCSIR$ from \eqref{eq:hDSIR}]{
\label{fig:quadratic1_SIR_evals}
\includegraphics[width=0.3\textwidth]{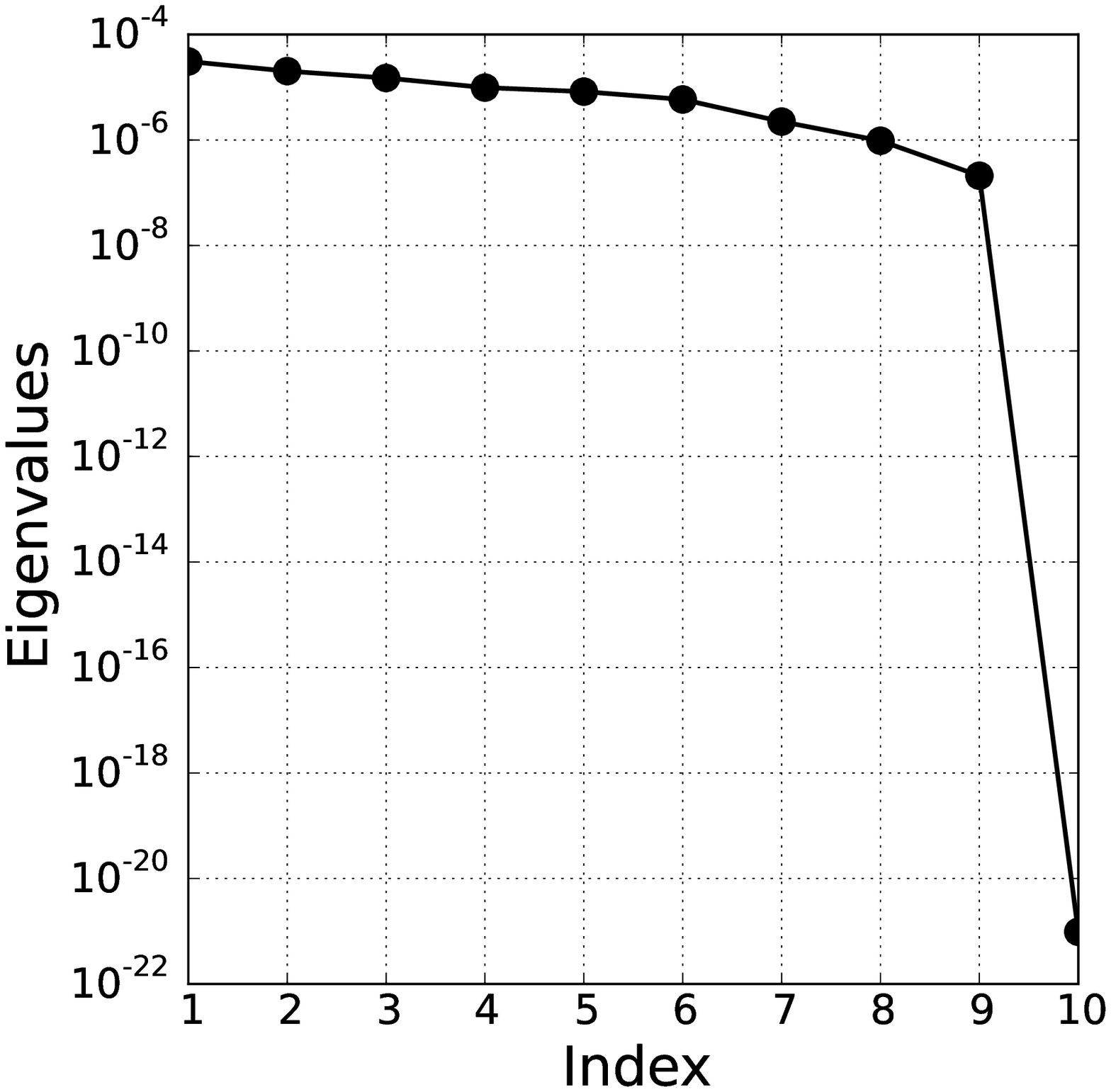}
}
\hfil
\subfloat[Sufficient summary plot for SIR]{
\label{fig:quadratic1_SIR_1d_ssp}
\includegraphics[width=0.3\textwidth]{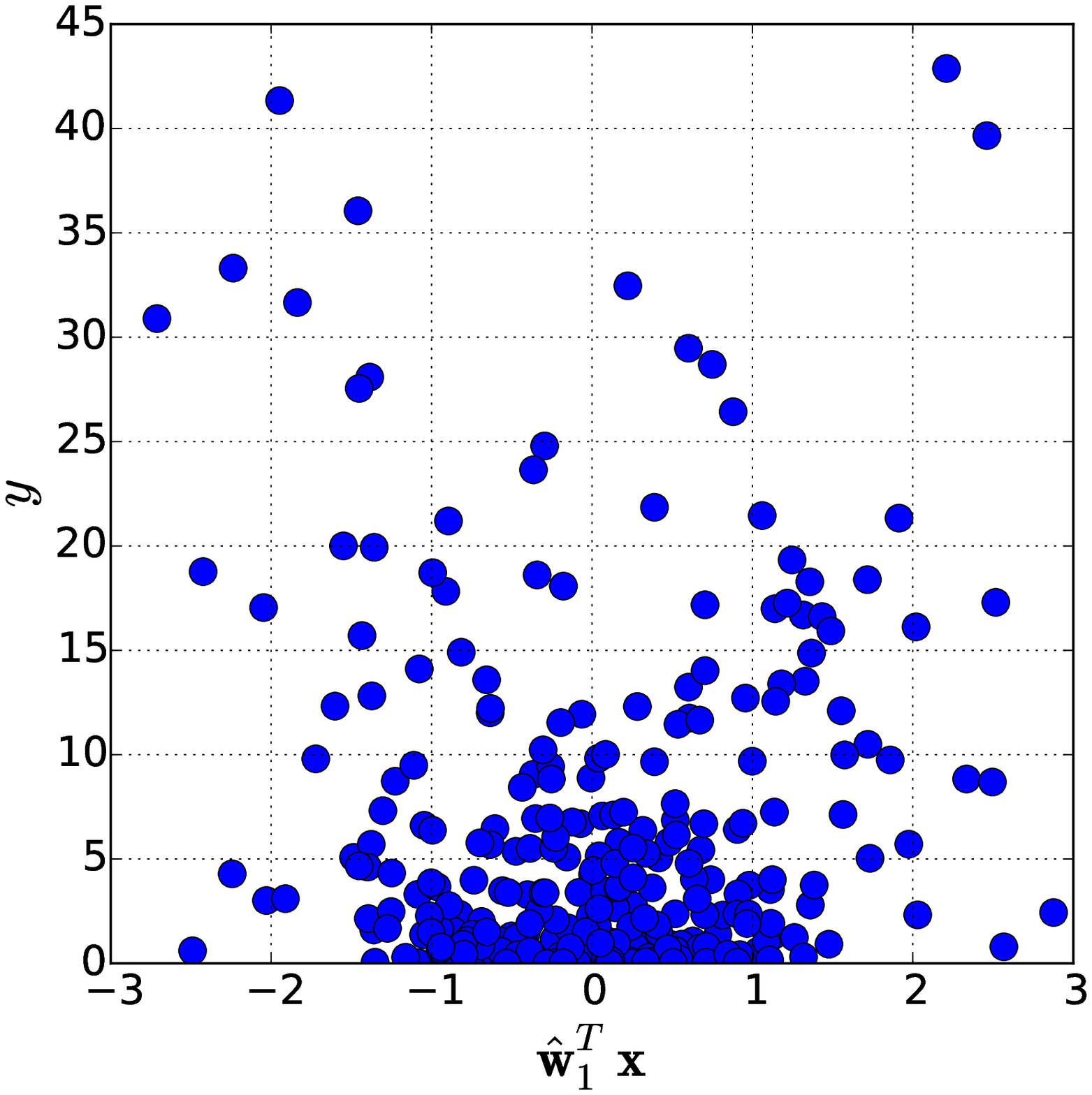}
}
\caption{As expected, SIR fails to recover the ridge direction $\vb$ in the function \eqref{eq:1d_quad}.}
\label{fig:quadratic1_SIR}

\end{figure*}

We study a simple one-dimensional quadratic ridge function to contrast the recovery properties of SIR versus SAVE. Let $\pi_{\vx}$ have a standard multivariate Gaussian density function on $\mathbb{R}^{10}$. Define
\begin{equation} \label{eq:1d_quad}
y \;=\; f(\vx) \;=\; \left( \vb^\top \vx \right)^2,
\end{equation}
where $\vb \in \mathbb{R}^{10}$ is a constant vector. The span of $\vb$ is the central subspace. First, we attempt to estimate the central subspace using SIR (Algorithm \ref{alg:SIR}), which is known to fail for functions symmetric about $\vx = \mzero$~\citep{Cook91}; Figure \ref{fig:quadratic1_SIR} confirms this failure. In fact, $\CIR$ is zero since the conditional expectation of $\vx$ for any value of $y$ is zero. Figure \ref{fig:quadratic1_SIR_evals} shows that all estimated eigenvalues of the SIR matrix are nearly zero as expected. Figure \ref{fig:quadratic1_SIR_1d_ssp} is a one-dimensional sufficient summary plot of $y_i$ against $\hat{\vw}_1^\top \vx_i$, where $\hat{\vw}_1$ denotes the normalized eigenvector associated with the largest eigenvalue of $\hCSIR$ from Algorithm \ref{alg:SIR}. If (i) the central subspace is one-dimensional (as in this case) and (ii) the chosen SDR algorithm correctly identifies the one basis vector, then the sufficient summary plot will show a univariate relationship between the linear combination of input evaluations and the associated outputs. Due to the symmetry in the quadratic function, SIR fails to recover the basis vector; the sufficient summary plot's lack of univariate relationship confirms the failure. 

\begin{figure*}[!ht]
\centering
\subfloat[Eigenvalues of $\hCSAVE$ from \eqref{eq:hDSAVE}]{
\label{fig:quadratic1_SAVE_evals}
\includegraphics[width=0.3\textwidth]{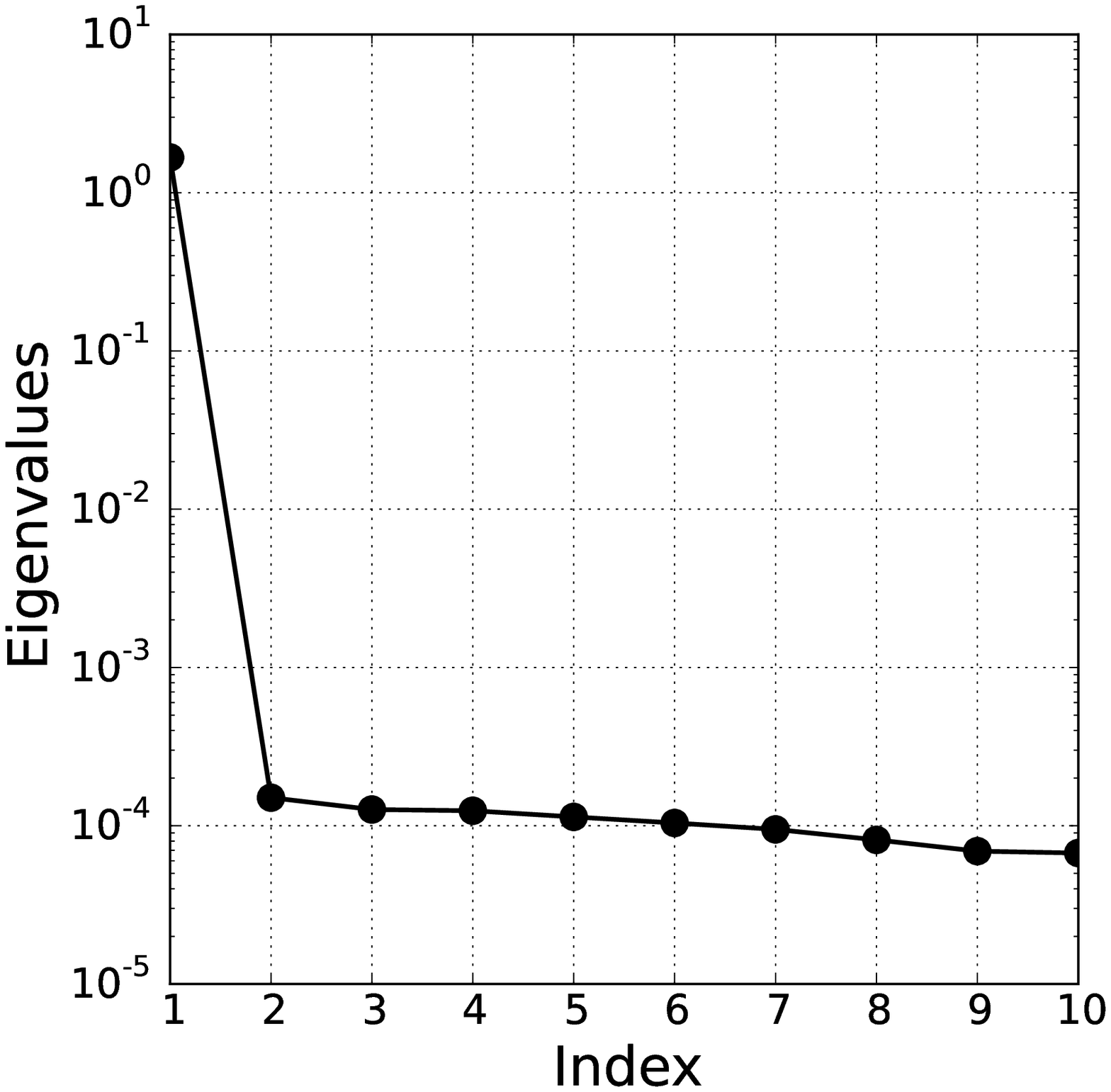}
}
\hfil
\subfloat[Sufficient summary plot for SAVE]{
\label{fig:quadratic1_SAVE_1d_ssp}
\includegraphics[width=0.3\textwidth]{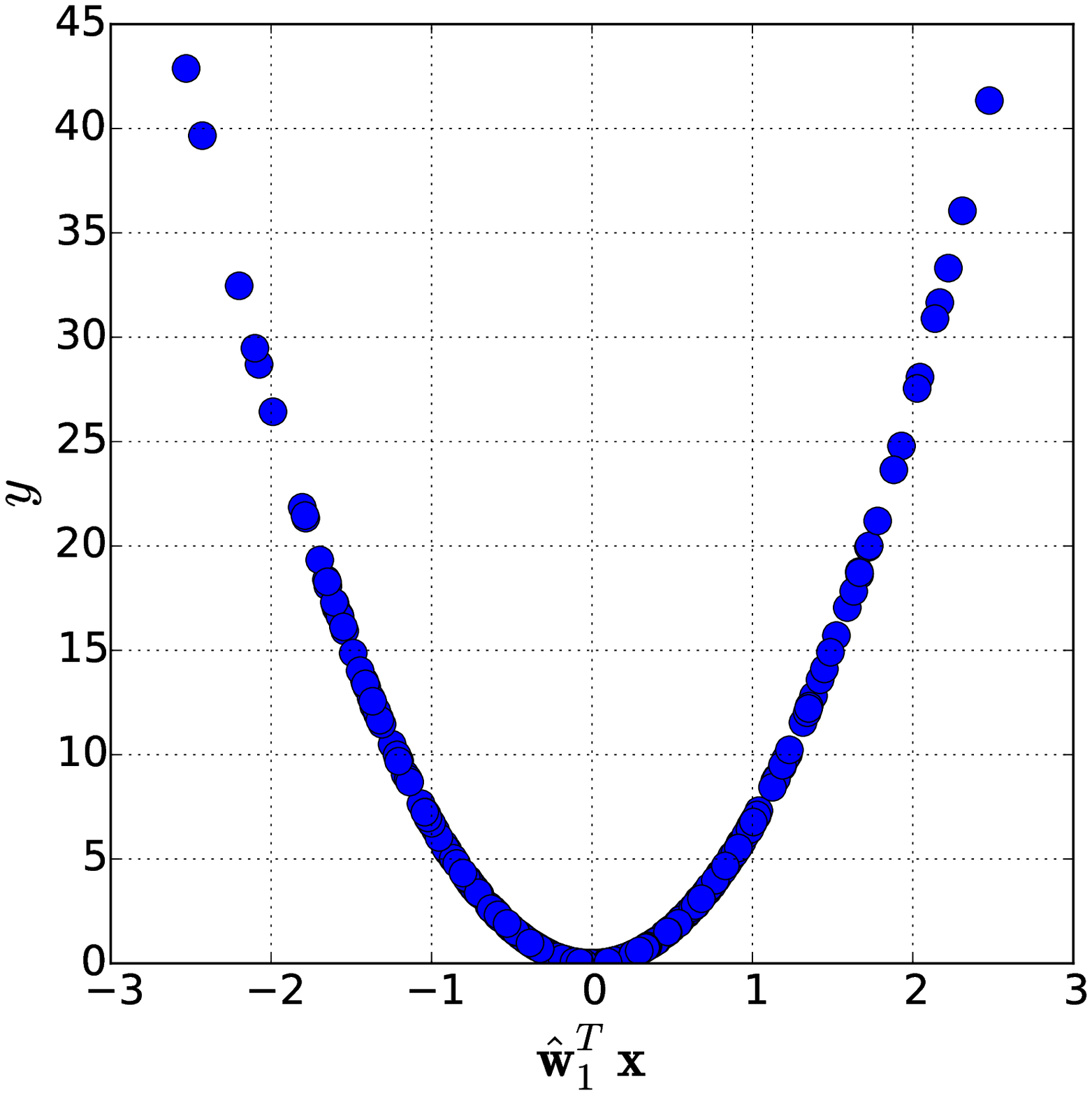}
}
\caption{SAVE recovers the ridge direction $\vb$ in the function \eqref{eq:1d_quad}.}
\label{fig:quadratic1_SAVE}
\end{figure*}

Figure \ref{fig:quadratic1_SAVE} shows results from applying SAVE (Algorithm \ref{alg:SAVE}) to the quadratic function \eqref{eq:1d_quad}. Figure \ref{fig:quadratic1_SAVE_evals} shows the eigenvalues of $\hCSAVE$ from Algorithm \ref{alg:SAVE}. Note the large gap between the first and second eigenvalues, which suggests that the SAVE subspace is one-dimensional. Figure \ref{fig:quadratic1_SAVE_1d_ssp} shows the sufficient summary plot using the first eigenvector $\hvw_1$ from Algorithm \ref{alg:SAVE}, which reveals the univariate quadratic relationship between $\hvw_1^\top\vx$ and $y$.

\subsection{Three-dimensional quadratic function}

\begin{figure*}[!ht]
\centering
\subfloat[Eigenvalues of $\hCSIR$ from \eqref{eq:hDSIR}]{
\label{fig:quadratic2_SIR_evals}
\includegraphics[width=0.3\textwidth]{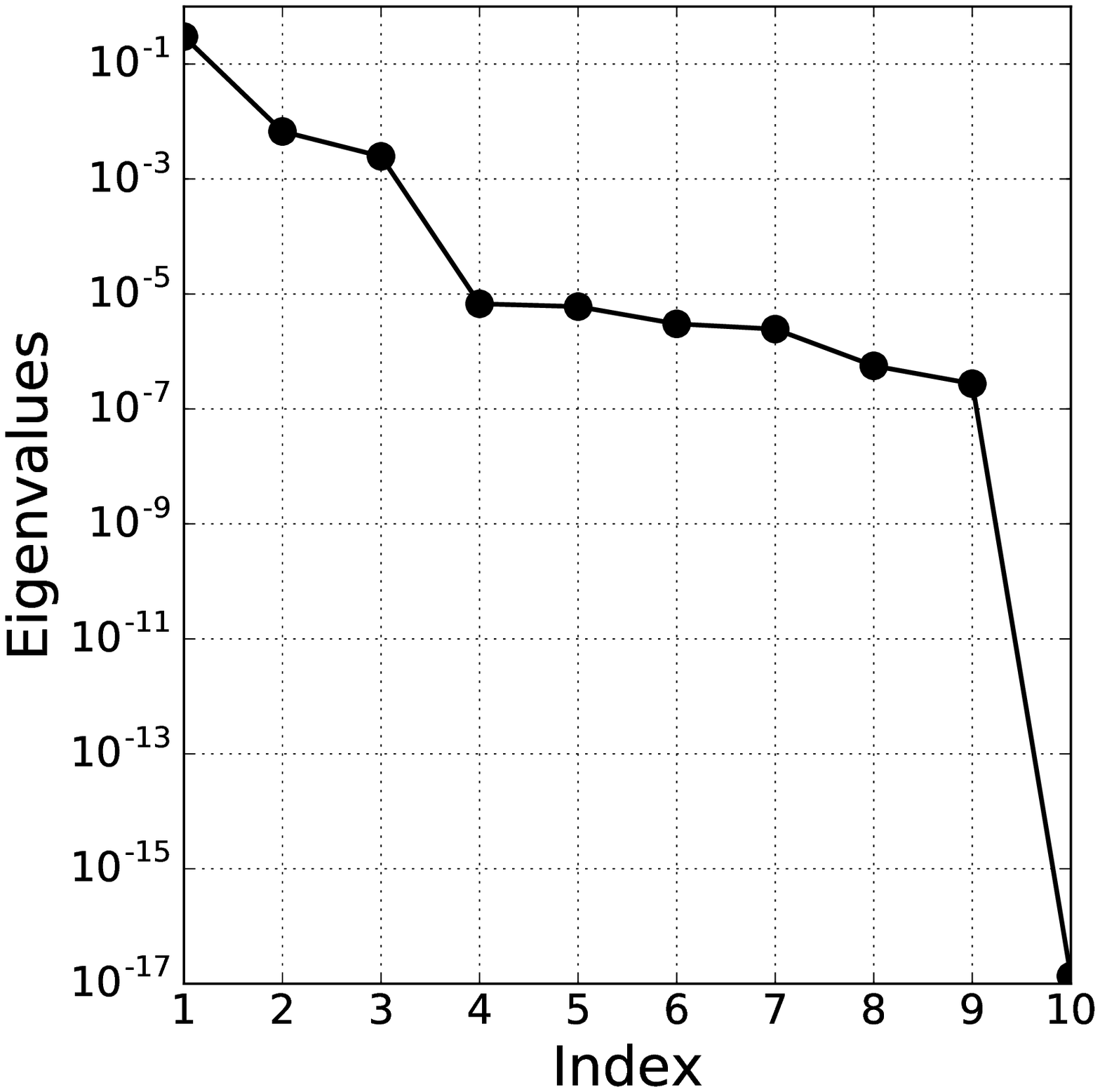}
}
\hfil
\subfloat[Maximum squared eigenvalue error, SIR]{
\label{fig:quadratic2_SIR_eval_errs}
\includegraphics[width=0.3\textwidth]{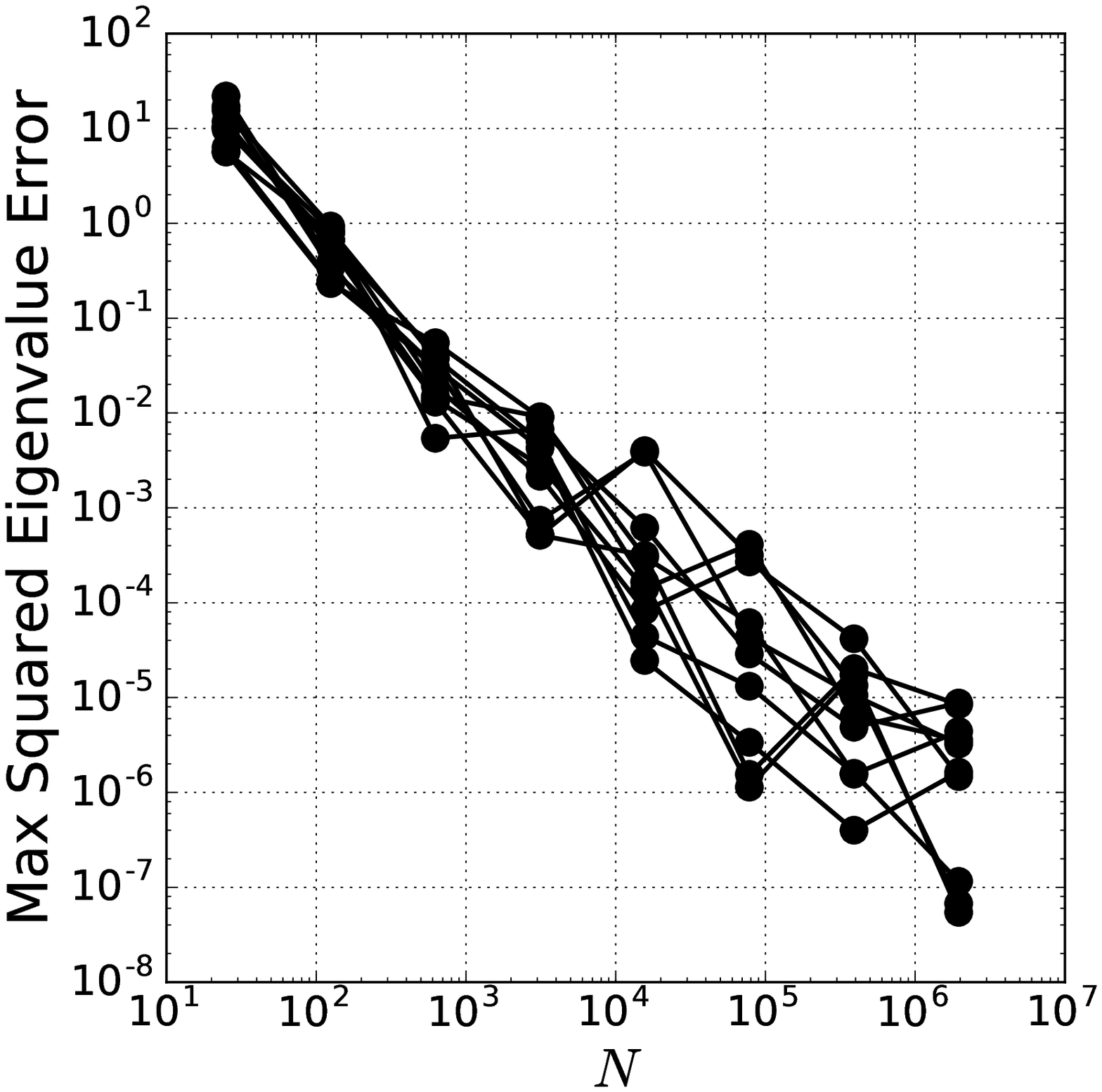}
}
\hfil
\subfloat[SIR subspace errors for $n = 3$]{
\label{fig:quadratic2_SIR_sub_errs}
\includegraphics[width=0.3\textwidth]{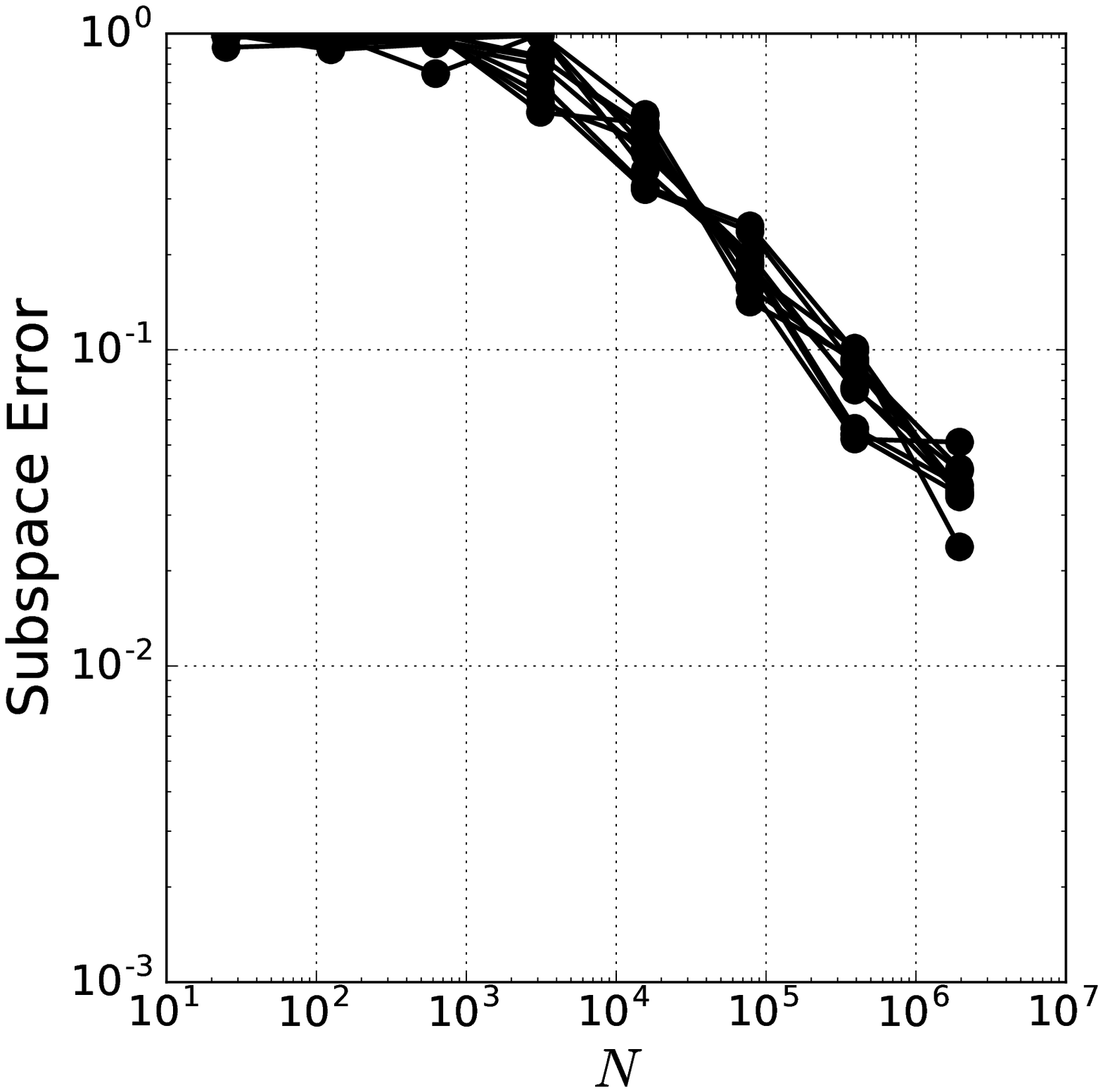}
}
\caption{Eigenvalues, eigenvalue errors, and subspace errors for SIR applied to \eqref{eq:3d_quad}. The error decreases with increasing samples consistent with the convergence theory in Section \ref{subsec:SIR_for_RR}.}
\label{fig:quadratic2_SIR}
\end{figure*}

Next, we numerically study the convergence properties of the SIR and SAVE algorithms using a more complex quadratic function. Let $\pi_{\vx}$ have a standard multivariate Gaussian density function on $\mathbb{R}^{10}$. Define
\begin{equation} \label{eq:3d_quad}
y \;=\; f(\vx) \;=\; \vx^\top \mB \mB^\top \vx + \vb^\top \vx,
\end{equation}
where $\mB \in \mathbb{R}^{10 \times 2}$ and $\vb \in \mathbb{R}^{10}$ with $\vb\not\in\colspan{\mB}$. Figure \ref{fig:quadratic2_SIR_evals} shows the eigenvalues of $\hCSIR$; note the gap between the third and fourth eigenvalues. Figure \ref{fig:quadratic2_SIR_eval_errs} shows the maximum squared eigenvalue error normalized by the largest eigenvalue,
\begin{equation}
\max_{1 \leq i \leq m} \frac{\left( \lambda_i (\hCSIR) - \lambda_i (\CSIR) \right)^2}{\lambda_1 (\CSIR)^2} ,
\end{equation}
for increasing numbers of samples in 10 independent trials. We estimate the true eigenvalues using SIR with $10^7$ samples. The average error decays at a rate slightly faster than the $\sO( N^{-1} )$ from Theorem \ref{thm:SIR_eig_converge}. The improvement can likely be attributed to the adaptive slicing procedure discussed at the beginning of this section. Figure \ref{fig:quadratic2_SIR_sub_errs} shows the error in the estimated three-dimensional SIR subspace (see \eqref{eq:sub_dist}) for increasing numbers of samples. We use $10^7$ samples to estimate the true SIR subspace. The subspace errors decrease asymptotically at a rate of approximately $\sO( N^{-1/2} )$, which agrees with Theorem \ref{thm:SIR_sub_converge}.

\begin{figure*}[!ht]
\centering
\subfloat[Eigenvalues of $\hCSAVE$ from \eqref{eq:hDSAVE}]{
\label{fig:quadratic2_SAVE_evals}
\includegraphics[width=0.25\textwidth]{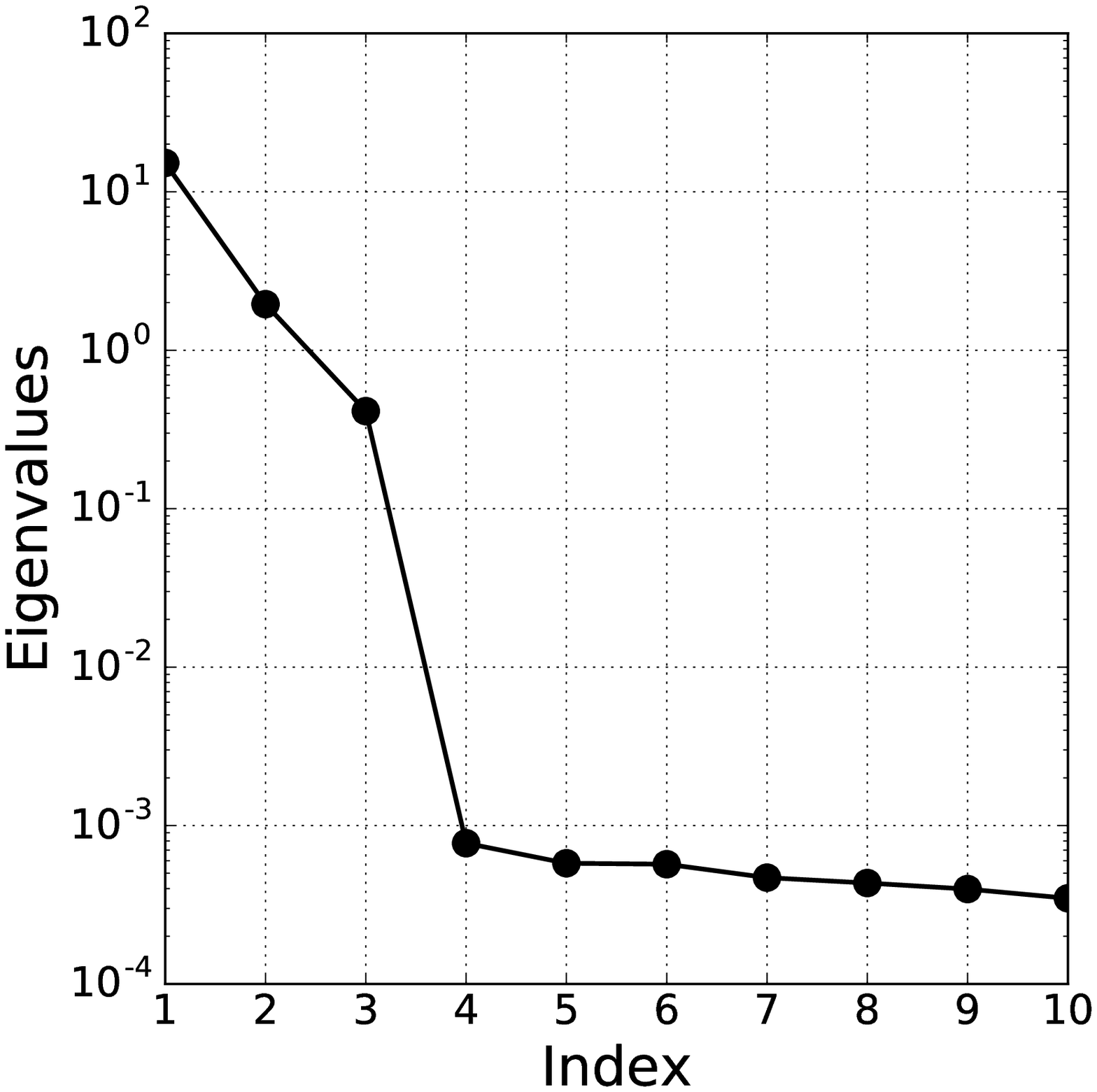}
} \hfil
\subfloat[Maximum squared eigenvalue error, SAVE]{
\label{fig:quadratic2_SAVE_eval_errs}
\includegraphics[width=0.25\textwidth]{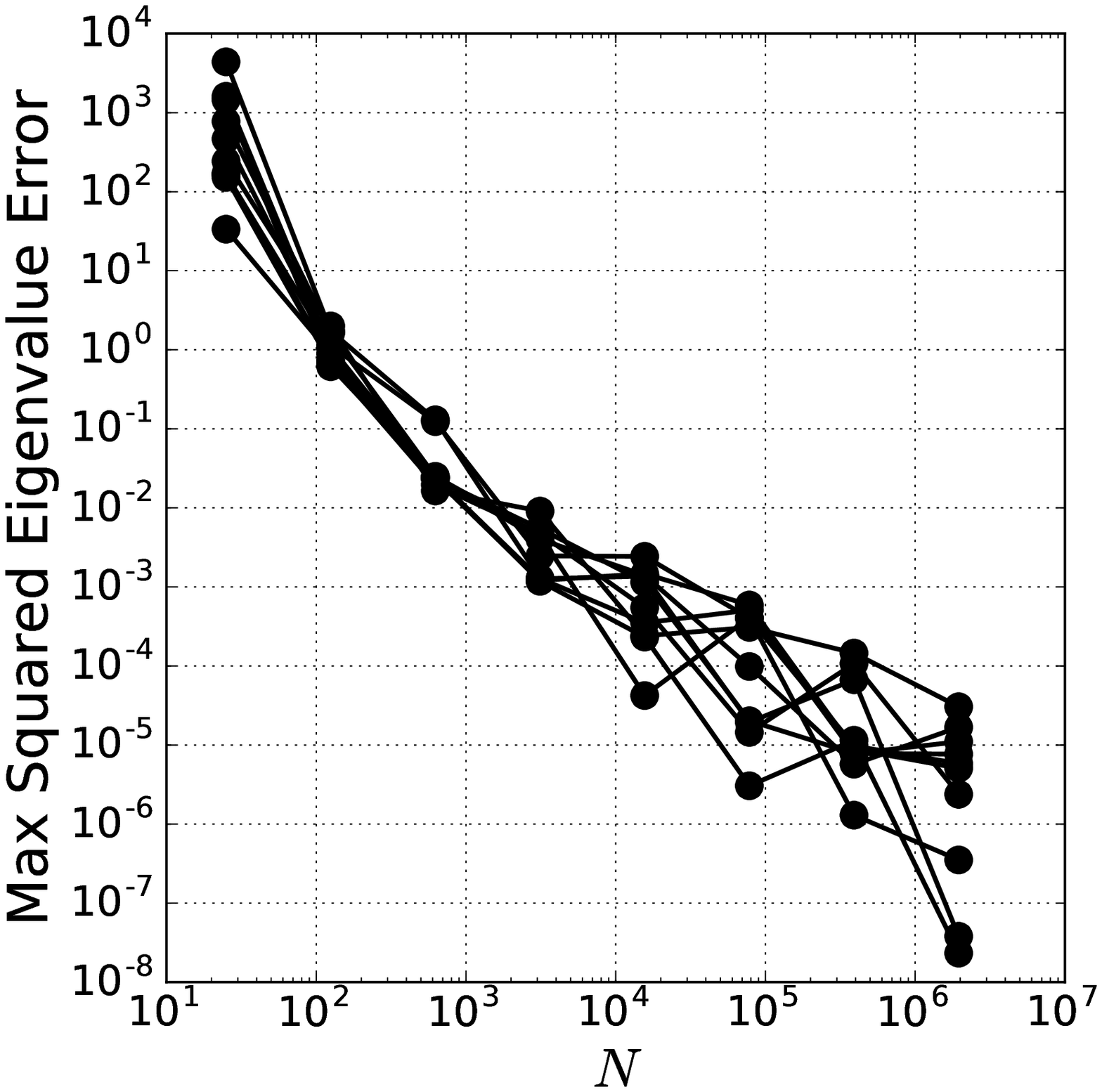}
}
\hfil
\subfloat[SAVE subspace errors for $n = 3$]{
\label{fig:quadratic2_SAVE_sub_errs}
\includegraphics[width=0.25\textwidth]{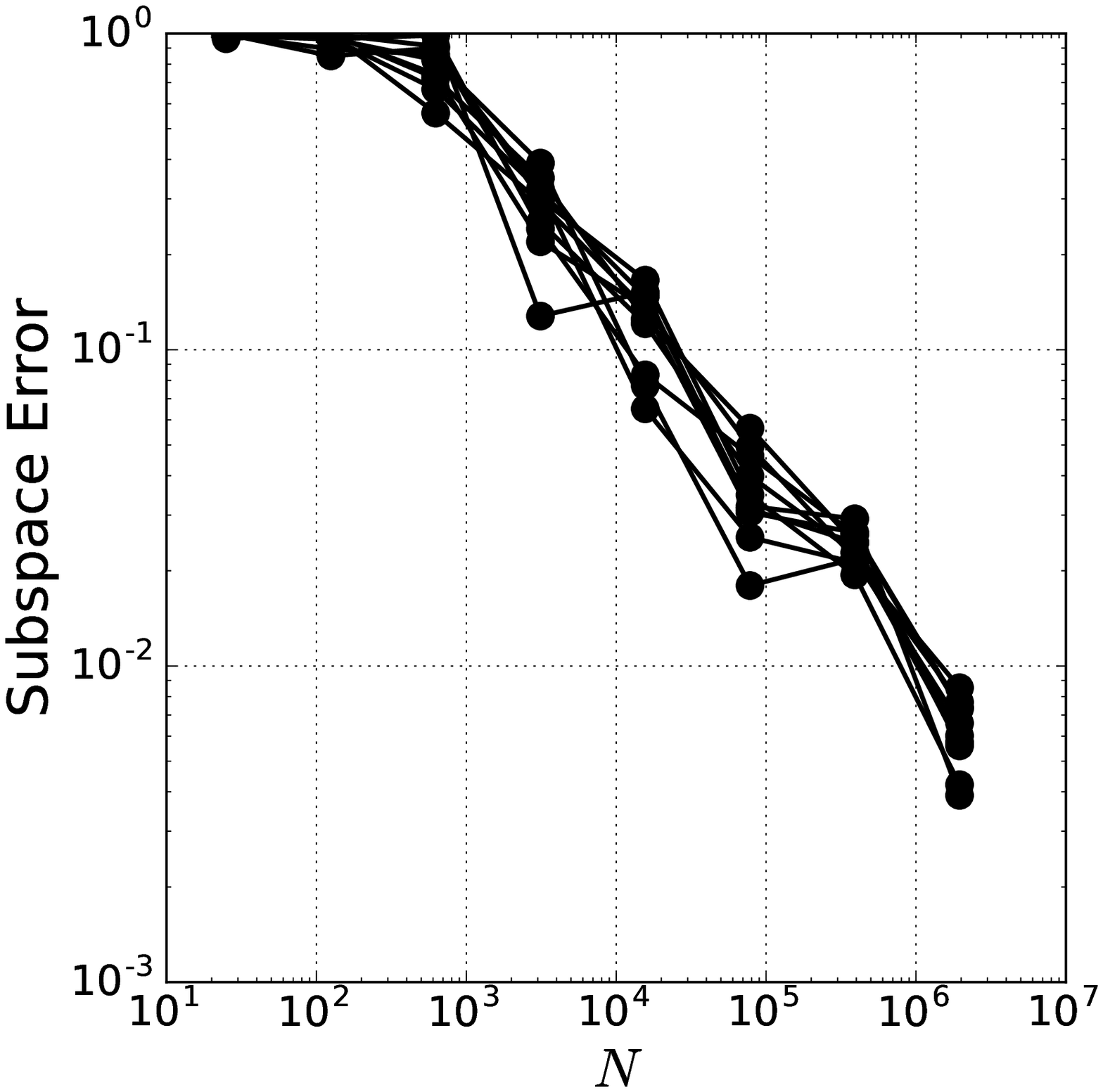}
}
\caption{Eigenvalues, eigenvalue errors, and subspace errors for SAVE applied to \eqref{eq:3d_quad}. The error decreases with increasing samples consistent with the convergence theory in Section \ref{subsec:SAVE_for_RR}.}
\label{fig:quadratic2_SAVE}
\end{figure*}

Figure \ref{fig:quadratic2_SAVE} shows the results of a similar convergence study using SAVE (Algorithm \ref{alg:SAVE}). The eigenvalues of $\hCSAVE$ from \eqref{eq:hDSAVE} are shown in Figure \ref{fig:quadratic2_SAVE_evals}. Note the large gap between the third and fourth eigenvalues, which is consistent with the three-dimensional central subspace in $f$ from \eqref{eq:3d_quad}. Figures \ref{fig:quadratic2_SAVE_eval_errs} and \ref{fig:quadratic2_SAVE_sub_errs} show the maximum squared eigenvalue error and the subspace error, respectively, for $n = 3$. The eigenvalue error again decays at a faster rate than expected in Theorem \ref{thm:SAVE_eig_converge}---likely due to the adaptive slicing implemented in the code. The subspace error decays consistently according to Theorem \ref{thm:SAVE_sub_converge}.

\subsection{Hartmann problem}
\label{subsec:hartmann}

\begin{figure}[ht]
\begin{center}
\includegraphics[width=0.4\textwidth]{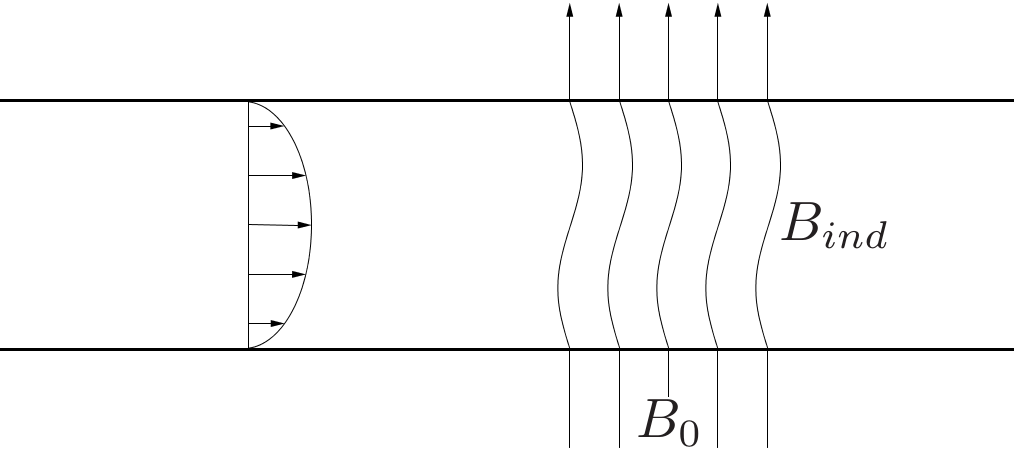}
\end{center}
\caption{The Hartmann problem studies the flow of an ionized fluid between two parallel plates. A magnetic field applied perpendicular to the flow direction acts as a resistive force to the fluid flow. Simultaneously, the fluid induces a magnetic field along the direction of the flow.}
\label{fig:hartmann_prob}
\end{figure}

The following study moves toward the direction of using SIR and SAVE for parameter reduction in a physics-based model.  The Hartmann problem is a standard problem in magnetohydrodynamics (MHD) that models the flow of an electrically-charged plasma in the presence of a uniform magnetic field~\citep{Cowling57}. The flow occurs along an infinite channel between two parallel plates separated by distance $2 \ell$. The applied magnetic field is perpendicular to the flow direction and acts as a resistive force on the flow velocity. At the same time, the movement of the fluid induces a magnetic field along the direction of the flow. Figure \ref{fig:hartmann_prob} contains a diagram of this problem. 
We have recently used the following model as a test case for parameter reduction methods~\citep{Glaws17}. 

\begin{figure*}[!ht]
\centering
\subfloat[Eigenvalues of $\hCSIR$ with bootstrap ranges]{
\label{fig:hartmann_Bind_SIR_evals}
\includegraphics[width=0.3\textwidth]{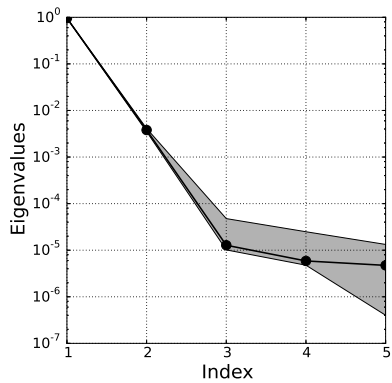}
}
\hfil
\subfloat[Subspace errors from $\hCSIR$]{
\label{fig:hartmann_Bind_SIR_suberr}
\includegraphics[width=0.3\textwidth]{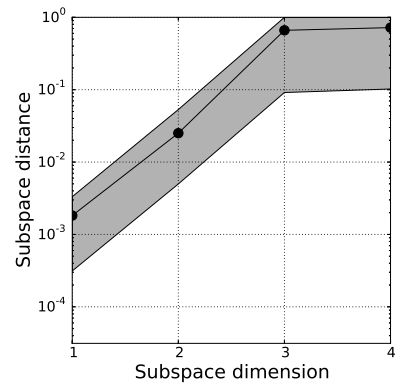}
} \\
\subfloat[One-dimensional SIR summary plot for $\Bind$ from \eqref{eq:Bind}]{
\label{fig:hartmann_Bind_SIR_1d_ssp}
\includegraphics[width=0.3\textwidth]{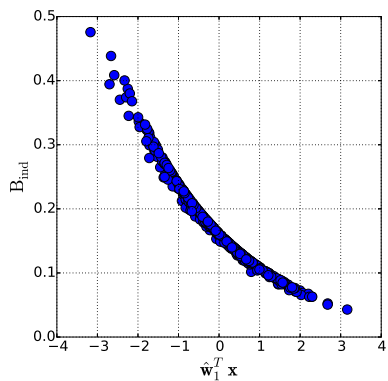}
}
\hfil
\subfloat[Two-dimensional SIR summary plot for $\Bind$ from \eqref{eq:Bind}]{
\label{fig:hartmann_Bind_SIR_2d_ssp}
\includegraphics[width=0.3\textwidth]{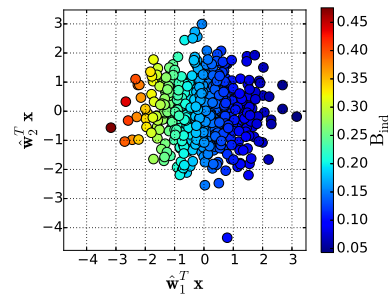}
}
\caption{Eigenvalues with bootstrap ranges, estimated subspace errors, and sufficient summary plots for SIR (Algorithm \ref{alg:SIR}) applied to $\Bind$ from \eqref{eq:Bind}.}
\label{fig:hartmann_Bind_SIR}
\end{figure*}

The inputs to the Hartmann model are fluid viscosity $\mu$, fluid density $\rho$, applied pressure gradient $\partial p_0 / \partial x$ (where the derivative is with respect to the flow field's spatial coordinate), resistivity $\eta$, and applied magnetic field $B_0$. We collect these inputs into a vector,
\begin{equation} \label{eq:hart_inputs}
\vx \;=\; \bmat{ \mu & \rho & \frac{\partial p_0}{\partial x} & \eta & B_0 }^\top.
\end{equation}
The output of interest is the total induced magnetic field,
\begin{equation}
\label{eq:Bind}
\Bind (\vx) = \frac{\partial p_0}{\partial x} \frac{\ell \mu_0}{2 B_0} \left( 1 - 2 \frac{\sqrt{\eta \mu}}{B_0 \ell} \, \text{tanh} \left( \frac{B_0 \ell}{2 \sqrt{\eta \mu}} \right) \right) .
\end{equation}
This function is not a ridge function of $\vx$. However, it has been shown that many physical laws can be expressed as ridge functions by considering a log transform of the inputs, which relates ridge structure in the function to dimension reduction via the Buckingham Pi Theorem~\citep{Constantine16}. For this reason, we apply SIR and SAVE as ridge recovery methods for $\Bind$ as a function of the logarithms of the inputs from \eqref{eq:hart_inputs}. The log-transformed inputs are equipped with a multivariate Gaussian with mean $\vmu$ and covariance $\mSigma$,
\begin{equation}
\vmu = \bmat{ -2.25 \\ 1 \\ 0.3 \\ 0.3 \\ -0.75 }, \qquad 
\mSigma = \bmat{ 0.15 & & & & \\ & 0.25 & & & \\ & & 0.25 & & \\ & & & 0.25 & \\ & & & & 0.25} .
\end{equation}
Figure \ref{fig:hartmann_Bind_SIR} shows the results of applying SIR (Algorithm \ref{alg:SIR}) to the Hartmann model for the induced magnetic field $\Bind$. The eigenvalues of $\hCSIR$ from \eqref{eq:hDSIR} with bootstrap ranges are shown in Figure \ref{fig:hartmann_Bind_SIR_evals}. Large gaps appear after the first and second eigenvalues, which indicates possible two-dimensional ridge structure. In fact, the induced magnetic field admits a two-dimensional central subspace relative to the log-inputs~\citep{Glaws17}. Figures \ref{fig:hartmann_Bind_SIR_1d_ssp} and \ref{fig:hartmann_Bind_SIR_2d_ssp} contain one- and two-dimensional sufficient summary plots of $\Bind$ against $\hvw_1^\top \vx$ and $\hvw_2^\top \vx$, where $\hvw_1$ and $\hvw_2$ are the first two eigenvectors of $\hCSIR$. We see a strong one-dimensional relationship. However, the two-dimensional sufficient summary plot shows slight curvature with changes in $\hvw_2^\top\vx$. These results suggest that ridge-like structure may be discovered using the SIR algorithm in some cases. Figure \ref{fig:hartmann_Bind_SIR_suberr} shows the subspace errors as a function of the subspace dimension. Recall from Theorem \ref{thm:SIR_sub_converge} that the subspace error depends inversely on the eigenvalue gap. The largest eigenvalue gap occurs between the first and second eigenvalues, which is consistent with the smallest subspace error for $n=1$.

\begin{figure*}[!ht]
\centering
\subfloat[Eigenvalues of $\hCSAVE$ with bootstrap ranges]{
\label{fig:hartmann_Bind_SAVE_evals}
\includegraphics[width=0.3\textwidth]{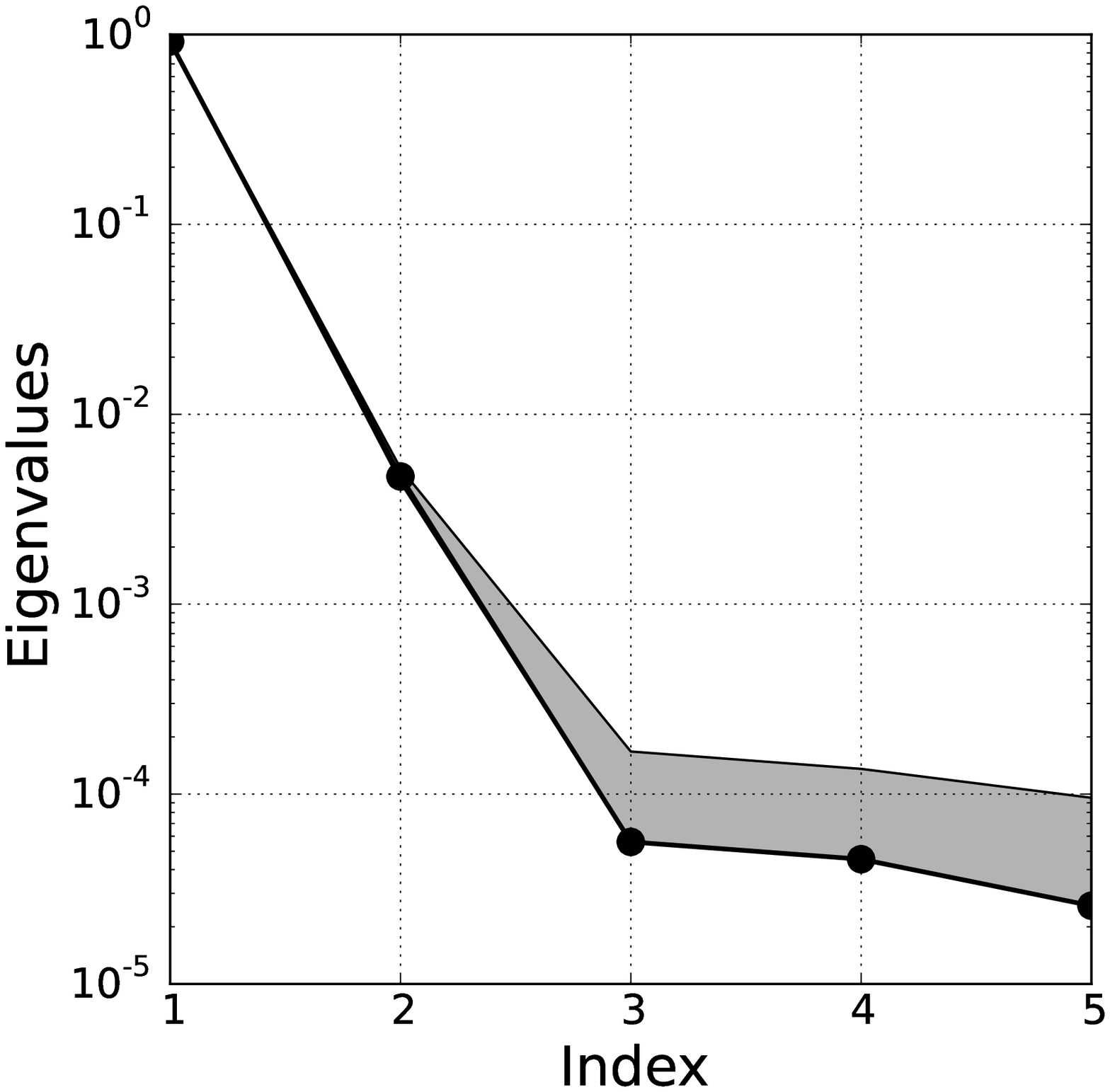}
}
\hfil
\subfloat[Subspace errors from $\hCSAVE$]{
\label{fig:hartmann_Bind_SAVE_suberr}
\includegraphics[width=0.3\textwidth]{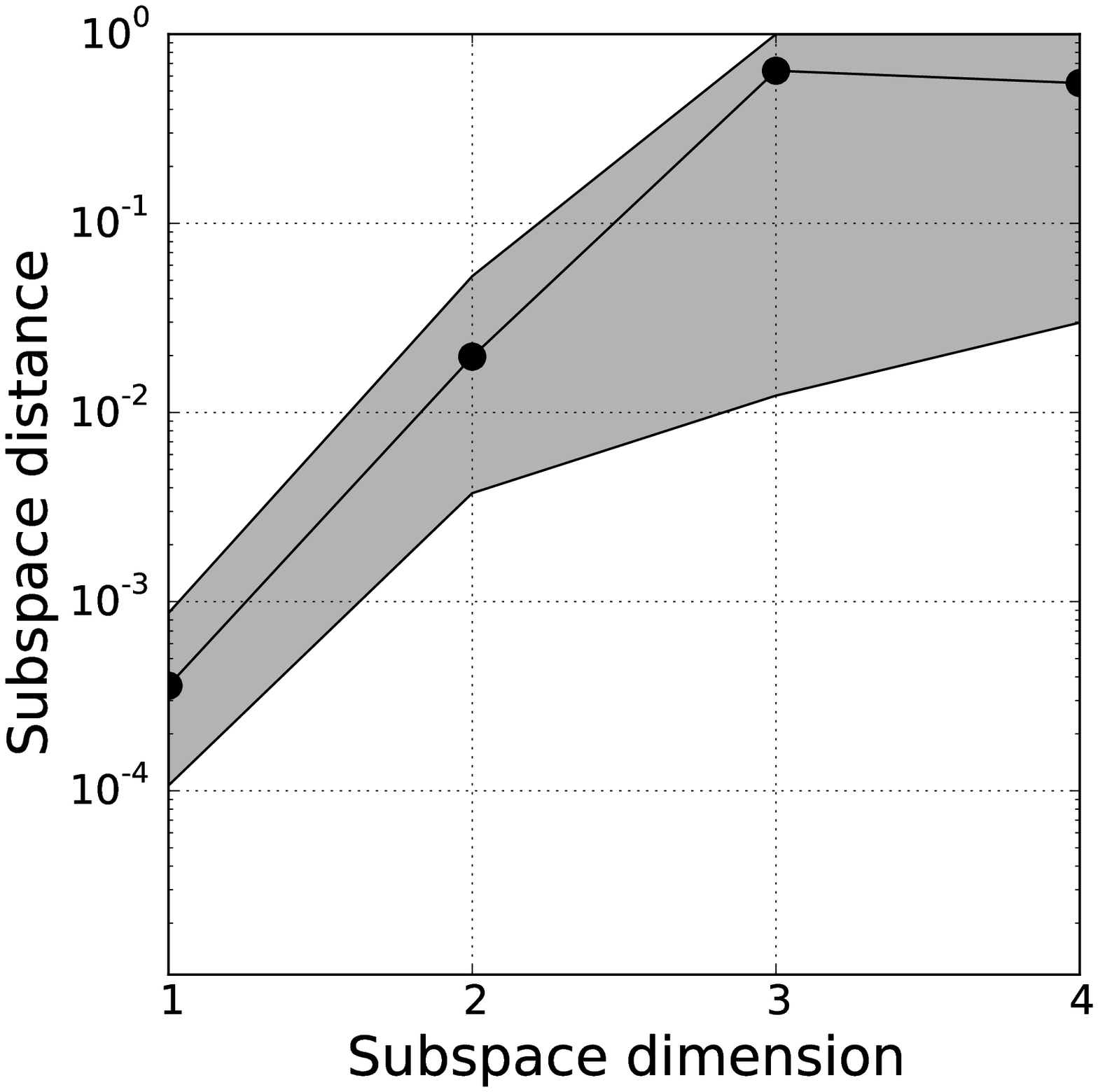}
} \\
\subfloat[One-dimensional SAVE summary plot for $\Bind$ from \eqref{eq:Bind}]{
\label{fig:hartmann_Bind_SAVE_1d_ssp}
\includegraphics[width=0.3\textwidth]{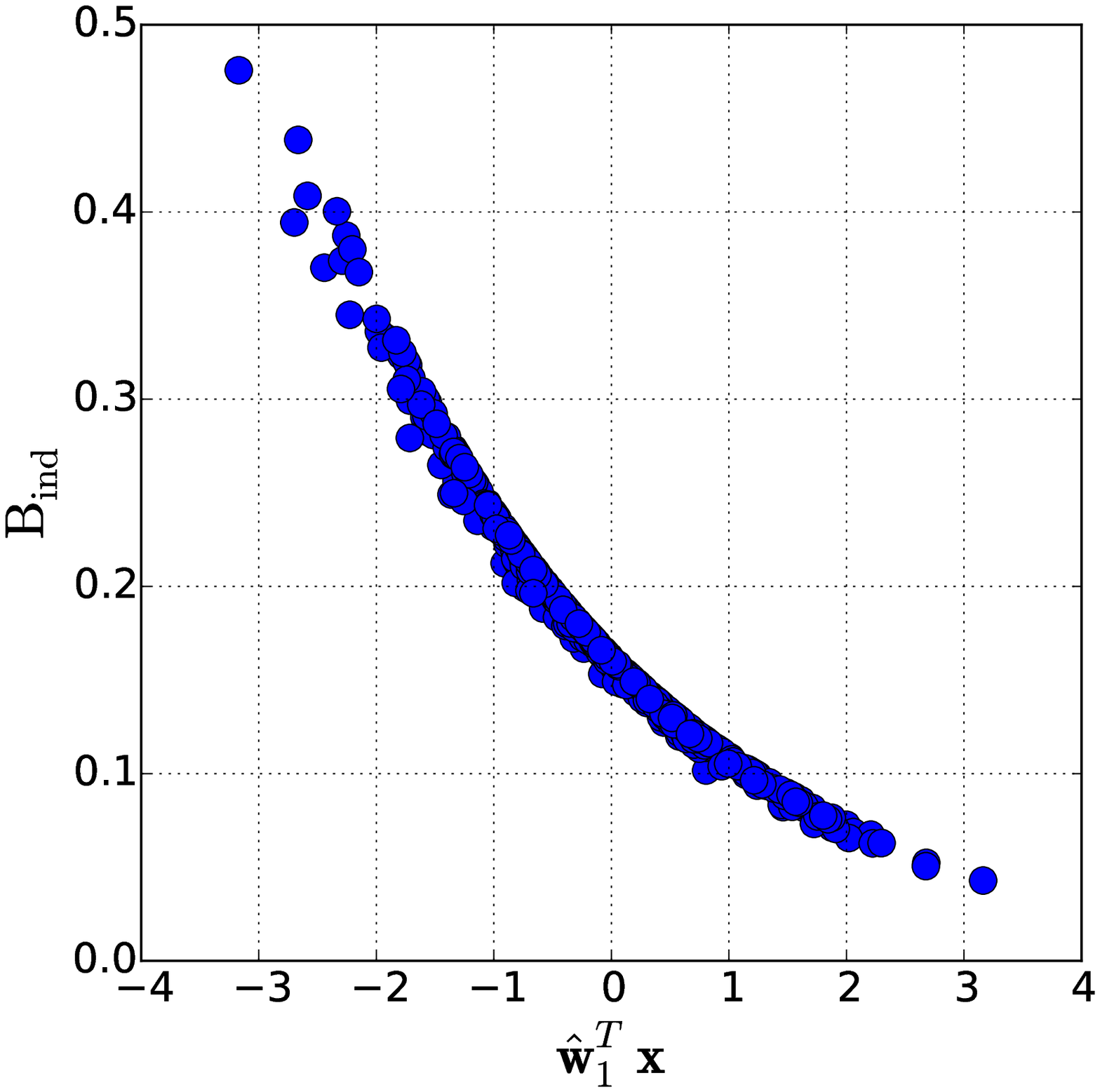}
}
\hfil
\subfloat[Two-dimensional SAVE summary plot for $\Bind$ from \eqref{eq:Bind}]{
\label{fig:hartmann_Bind_SAVE_2d_ssp}
\includegraphics[width=0.3\textwidth]{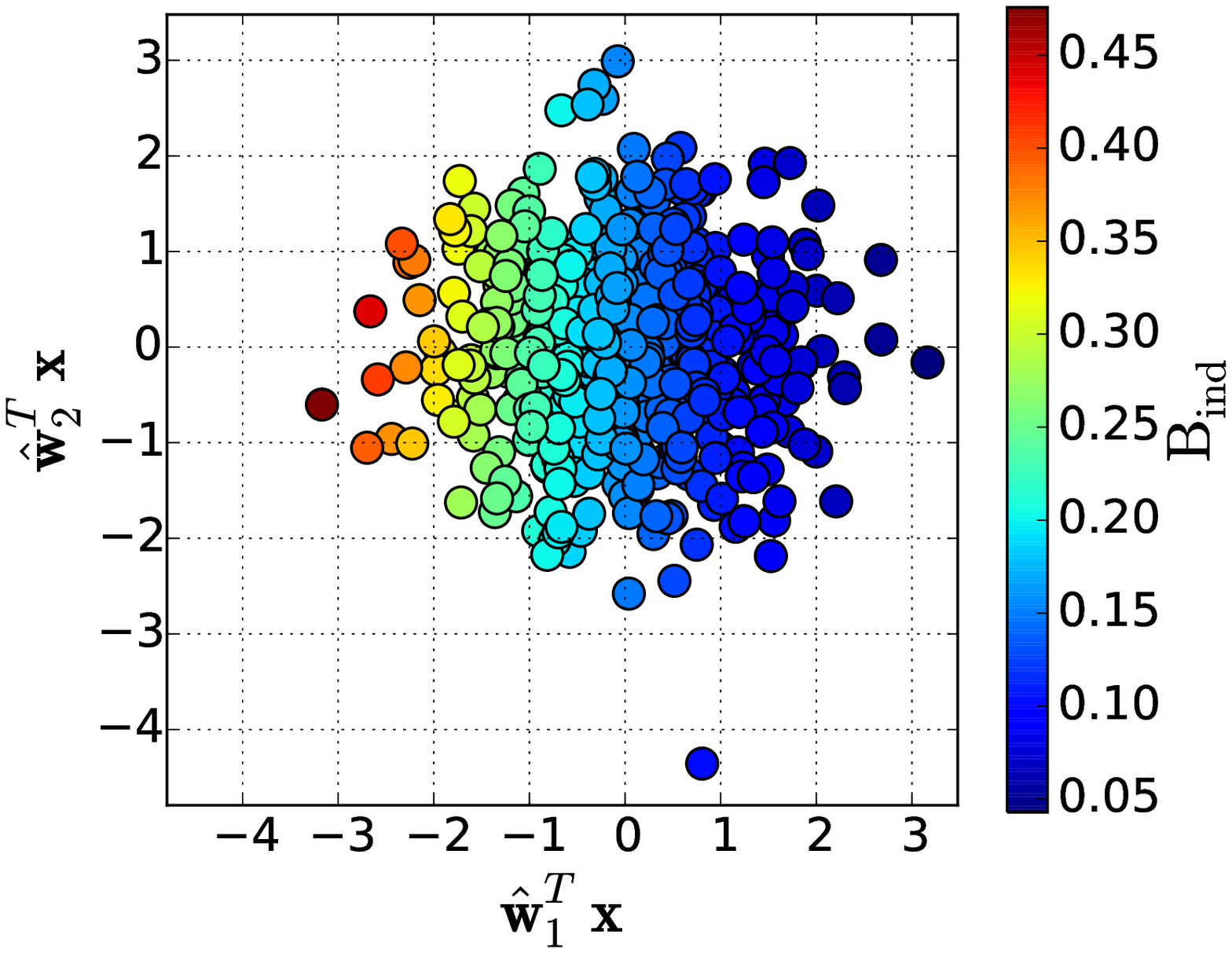}
}
\caption{Eigenvalues with bootstrap ranges, estimated subspace errors, and sufficient summary plots for SAVE (Algorithm \ref{alg:SAVE}) applied to $\Bind$ from \eqref{eq:Bind}.}
\label{fig:hartmann_Bind_SAVE}
\end{figure*}

We perform the same numerical studies using SAVE. Figure \ref{fig:hartmann_Bind_SAVE_evals} shows the eigenvalues of the $\hCSAVE$ from \eqref{eq:hDSAVE} for the induced magnetic field $\Bind$ from \eqref{eq:Bind}. Note the large gaps after the first and second eigenvalues. These gaps are consistent with the subspace errors in Figure \ref{fig:hartmann_Bind_SAVE_suberr}, where the one- and two-dimensional subspace estimates have the smallest errors. Figures \ref{fig:hartmann_Bind_SAVE_1d_ssp} and \ref{fig:hartmann_Bind_SAVE_2d_ssp} contain sufficient summary plots for $\hvw_1^\top\vx$ and $\hvw_2^\top\vx$, where $\hvw_1$ and $\hvw_2$ are the first two eigenvectors from $\hCSAVE$ in \eqref{eq:hDSAVE}.

\section{Summary and conclusion}
\label{sec:conclusion}

We investigate sufficient dimension reduction from statistical regression as a tool for subspace-based parameter reduction in deterministic functions where applications of interest include computer experiments. We show that SDR is theoretically justified as a tool for ridge recovery by proving equivalence of the dimension reduction subspace and the ridge subspace for some deterministic $y = f(\vx)$. We interpret two SDR algorithms for the ridge recovery problem: sliced inverse regression and sliced average variance estimation. In regression, these methods use moments of the inverse regression $\vx|y$ to estimate subspaces relating to the central subspace. In ridge recovery, we reinterpret SIR and SAVE as numerical integration methods for estimating inverse conditional moment matrices, where the integrals are over contour sets of $f$. We show that the column spaces of the conditional moment matrices are contained in the ridge subspace, which justifies their eigenspaces as tools ridge recovery.

\bibliographystyle{natbib}
\bibliography{deterministicSDR}

\begin{thebibliography}{}

\bibitem[Abbott(2001)Abbott]{Abbott01}
Abbott, S. (2001).
\newblock {\em Understanding Analysis\/}.
\newblock Springer, New York, 2nd edition.

\bibitem[Adragni and Cook(2009)Adragni and Cook]{Adragni09}
Adragni, K.~P. and Cook, R.~D. (2009).
\newblock Sufficient dimension reduction and prediction in regression.
\newblock {\em Philosophical Transactions of the Royal Society of London A:
  Mathematical, Physical and Engineering Sciences\/}, {\bf 367}(1906),
  4385--4405.

\bibitem[Allaire and Willcox(2010)Allaire and Willcox]{Allaire2010}
Allaire, D. and Willcox, K. (2010).
\newblock Surrogate modeling for uncertainty assessment with application to
  aviation environmental system models.
\newblock {\em AIAA Journal\/}, {\bf 48}(8), 1791--1803.

\bibitem[Challenor(2012)Challenor]{Challenor2012}
Challenor, P. (2012).
\newblock Using emulators to estimate uncertainty in complex models.
\newblock In A.~M. Dienstfrey and R.~F. Boisvert, editors, {\em Uncertainty
  Quantification in Scientific Computing\/}, pages 151--164, Berlin,
  Heidelberg. Springer Berlin Heidelberg.

\bibitem[Chang and Pollard(1997)Chang and Pollard]{Chang1997}
Chang, J.~T. and Pollard, D. (1997).
\newblock Conditioning as disintegration.
\newblock {\em Statistica Neerlandica\/}, {\bf 51}(3), 287--317.

\bibitem[Cohen {\em et~al.}(2012)Cohen, Daubechies, DeVore, Kerkyacharian, and
  Picard]{Cohen12}
Cohen, A., Daubechies, I., DeVore, R., Kerkyacharian, G., and Picard, D.
  (2012).
\newblock Capturing ridge functions in high dimensions from point queries.
\newblock {\em Constructive Approximation\/}, {\bf 35}(2), 225--243.

\bibitem[Constantine(2015)Constantine]{Constantine15}
Constantine, P.~G. (2015).
\newblock {\em Active Subspaces: Emerging Ideas for Dimension Reduction in
  Parameter Studies\/}.
\newblock SIAM, Philadelphia.

\bibitem[Constantine {\em et~al.}(2015)Constantine, Eftekhari, and
  Wakin]{constantine2015sketching}
Constantine, P.~G., Eftekhari, A., and Wakin, M.~B. (2015).
\newblock Computing active subspaces efficiently with gradient sketching.
\newblock In {\em 2015 IEEE 6th International Workshop on Computational
  Advances in Multi-Sensor Adaptive Processing (CAMSAP)\/}, pages 353--356.

\bibitem[Constantine {\em et~al.}(2016)Constantine, Howard, Glaws, Grey, Diaz,
  and Fletcher]{constantine2016python}
Constantine, P.~G., Howard, R., Glaws, A., Grey, Z., Diaz, P., and Fletcher, L.
  (2016).
\newblock Python active-subspaces utility library.
\newblock {\em The Journal of Open Source Software\/}, {\bf 1}.

\bibitem[Constantine {\em et~al.}(2017a)Constantine, del Rosario, and
  Iaccarino]{Constantine16}
Constantine, P.~G., del Rosario, Z., and Iaccarino, G. (2017a).
\newblock Data-driven dimensional analysis: algorithms for unique and relevant
  dimensionless groups.
\newblock {\em arXiv:1708.04303\/}.

\bibitem[Constantine {\em et~al.}(2017b)Constantine, Eftekhari, Hokanson, and
  Ward]{constantine2016near}
Constantine, P.~G., Eftekhari, A., Hokanson, J., and Ward, R. (2017b).
\newblock A near-stationary subspace for ridge approximation.
\newblock {\em Computer Methods in Applied Mechanics and Engineering\/}, {\bf
  326}, 402--421.

\bibitem[Cook(1994a)Cook]{Cook94a}
Cook, R.~D. (1994a).
\newblock On the interpretation of regression plots.
\newblock {\em Journal of the American Statistical Association\/}, {\bf
  89}(425), 177--189.

\bibitem[Cook(1994b)Cook]{Cook94b}
Cook, R.~D. (1994b).
\newblock Using dimension-reduction subspaces to identify important inputs in
  models of physical systems.
\newblock In {\em Proceedings of the Section on Physical and Engineering
  Sciences\/}, pages 18--25. American Statistical Association, Alexandria, VA.

\bibitem[Cook(1996)Cook]{Cook96}
Cook, R.~D. (1996).
\newblock Graphics for regressions with a binary response.
\newblock {\em Journal of the American Statistical Association\/}, {\bf
  91}(435), 983--992.

\bibitem[Cook(1998)Cook]{Cook98}
Cook, R.~D. (1998).
\newblock {\em Regression Graphics: Ideas for Studying Regression through
  Graphics\/}.
\newblock John Wiley \& Sons, Inc, New York.

\bibitem[Cook(2000)Cook]{Cook00a}
Cook, R.~D. (2000).
\newblock {SAVE}: A method for dimension reduction and graphics in regression.
\newblock {\em Communications in Statistics - Theory and Methods\/}, {\bf
  29}(9-10), 2109--2121.

\bibitem[Cook and Forzani(2009)Cook and Forzani]{Cook09}
Cook, R.~D. and Forzani, L. (2009).
\newblock Likelihood-based sufficient dimension reduction.
\newblock {\em Journal of the American Statistical Association\/}, {\bf
  104}(485), 197--208.

\bibitem[Cook and Weisberg(1991)Cook and Weisberg]{Cook91}
Cook, R.~D. and Weisberg, S. (1991).
\newblock Sliced inverse regression for dimension reduction: comment.
\newblock {\em Journal of the American Statistical Association\/}, {\bf
  86}(414), 328--332.

\bibitem[Cowling and Lindsay(1957)Cowling and Lindsay]{Cowling57}
Cowling, T. and Lindsay, R.~B. (1957).
\newblock Magnetohydrodynamics.
\newblock {\em Physics Today\/}, {\bf 10}, 40.

\bibitem[Donoho(2000)Donoho]{Donoho00}
Donoho, D.~L. (2000).
\newblock High-dimensional data analysis: The curses and blessings of
  dimensionality.
\newblock In {\em AMS Conference on Math Challenges of the 21st Century\/}.

\bibitem[Eaton(1986)Eaton]{Eaton86}
Eaton, M.~L. (1986).
\newblock A characterization of spherical distributions.
\newblock {\em Journal of Multivariate Analysis\/}, {\bf 20}(2), 272--27.

\bibitem[Folland(1999)Folland]{Folland99}
Folland, G.~B. (1999).
\newblock {\em Real Analysis: Modern Techniques and Their Applications\/}.
\newblock John Wiley \& Sons Ltd, New York, 2nd edition.

\bibitem[Fornasier {\em et~al.}(2012)Fornasier, Schnass, and
  Vybiral]{Fornaiser12}
Fornasier, M., Schnass, K., and Vybiral, J. (2012).
\newblock Learning functions of few arbitrary linear parameters in high
  dimensions.
\newblock {\em Foundations of Computational Mathematics\/}, {\bf 12}(2),
  229--262.

\bibitem[Friedman and Stuetzle(1980)Friedman and Stuetzle]{Friedman80}
Friedman, J.~H. and Stuetzle, W. (1980).
\newblock Projection pursuit regression.
\newblock {\em Journal of the American Statistical Association\/}, {\bf
  76}(376), 817--823.

\bibitem[Ghanem {\em et~al.}(2016)Ghanem, Higdon, and Owhadi]{HandbookUQ}
Ghanem, R., Higdon, D., and Owhadi, H. (2016).
\newblock {\em Handbook of Uncertainty Quantification\/}.
\newblock Springer International Publishing.

\bibitem[Glaws and Constantine(2018)Glaws and Constantine]{Glaws2018}
Glaws, A. and Constantine, P.~G. (2018).
\newblock {Gauss--Christoffel} quadrature for inverse regression: applications
  to computer experiments.
\newblock {\em Statistics and Computing\/}.

\bibitem[Glaws {\em et~al.}(2017)Glaws, Constantine, Shadid, and
  Wildey]{Glaws17}
Glaws, A., Constantine, P.~G., Shadid, J., and Wildey, T.~M. (2017).
\newblock Dimension reduction in magnetohydrodynamics power generation models:
  dimensional analysis and active subspaces.
\newblock {\em Statistical Analysis and Data Mining\/}, {\bf 10}(5), 312--325.

\bibitem[Golub and Van~Loan(2013)Golub and Van~Loan]{Golub96}
Golub, G.~H. and Van~Loan, C.~F. (2013).
\newblock {\em Matrix Computations\/}.
\newblock JHU Press, Baltimore, 4th edition.

\bibitem[Goodfellow {\em et~al.}(2016)Goodfellow, Bengio, and
  Courville]{Goodfellow16}
Goodfellow, I., Bengio, Y., and Courville, A. (2016).
\newblock {\em Deep Learning\/}.
\newblock MIT Press, Cambridge.

\bibitem[Hastie {\em et~al.}(2009)Hastie, Tibshirani, and Friedman]{Hastie2009}
Hastie, T., Tibshirani, R., and Friedman, J. (2009).
\newblock {\em The Elements of Statistical Learning: Data Mining, Inference,
  and Prediction\/}.
\newblock Springer, New York, 2nd edition.

\bibitem[Hokanson and Constantine(2017)Hokanson and Constantine]{Hokanson17}
Hokanson, J. and Constantine, P. (2017).
\newblock Data-driven polynomial ridge approximation using variable projection.
\newblock {\em arXiv:1702.05859\/}.

\bibitem[Jones(2001)Jones]{Jones2001}
Jones, D.~R. (2001).
\newblock A taxonomy of global optimization methods based on response surfaces.
\newblock {\em Journal of Global Optimization\/}, {\bf 21}(4), 345--383.

\bibitem[Koehler and Owen(1996)Koehler and Owen]{Koehler96}
Koehler, J.~R. and Owen, A.~B. (1996).
\newblock Computer experiments.
\newblock {\em Handbook of Statistics\/}, {\bf 13}(9), 261--308.

\bibitem[Li(2018)Li]{Li2018}
Li, B. (2018).
\newblock {\em Sufficient Dimension Reduction: Methods and Applications with
  R\/}.
\newblock CRC Press, Philadelphia.

\bibitem[Li(1991)Li]{Li91}
Li, K.~C. (1991).
\newblock Sliced inverse regression for dimension reduction.
\newblock {\em Journal of the American Statistical Association\/}, {\bf
  86}(414), 316--327.

\bibitem[Li(1992)Li]{Li92}
Li, K.~C. (1992).
\newblock On principal hessian directions for data visualization and dimension
  reduction: Another application of stein's lemma.
\newblock {\em Journal of the American Statistical Association\/}, {\bf
  87}(420), 1025--1039.

\bibitem[Li and Duan(1989)Li and Duan]{OLS89}
Li, K.-C. and Duan, N. (1989).
\newblock Regression analysis under link violation.
\newblock {\em The Annals of Statistics\/}, {\bf 17}(3), 1009--1052.

\bibitem[Li {\em et~al.}(2016)Li, Lin, and Li]{Li16}
Li, W., Lin, G., and Li, B. (2016).
\newblock Inverse regression-based uncertainty quantification algorithms for
  high-dimensional models: theory and practice.
\newblock {\em Journal of Computational Physics\/}, {\bf 321}, 259--278.

\bibitem[Myers and Montgomery(1995)Myers and Montgomery]{Myers1995}
Myers, R.~H. and Montgomery, D.~C. (1995).
\newblock {\em Response Surface Methodology: Process and Product Optimization
  Using Designed Experiments\/}.
\newblock John Wiley \& Sons, New York.

\bibitem[Pan and Dias(2017)Pan and Dias]{Pan2017}
Pan, Q. and Dias, D. (2017).
\newblock Sliced inverse regression-based sparse polynomial chaos expansions
  for reliability analysis in high dimensions.
\newblock {\em Reliability Engineering \& System Safety\/}, {\bf 167},
  484--493.
\newblock Special Section: Applications of Probabilistic Graphical Models in
  Dependability, Diagnosis and Prognosis.

\bibitem[Pinkus(2015)Pinkus]{Pinkus15}
Pinkus, A. (2015).
\newblock {\em Ridge Functions\/}.
\newblock Cambridge University Press.

\bibitem[Razavi {\em et~al.}(2012)Razavi, Tolson, and Burn]{Razavi2012}
Razavi, S., Tolson, B.~A., and Burn, D.~H. (2012).
\newblock Review of surrogate modeling in water resources.
\newblock {\em Water Resources Research\/}, {\bf 48}(7), W07401.

\bibitem[Sacks {\em et~al.}(1989)Sacks, Welch, Mitchell, and Wynn]{Sacks89}
Sacks, J., Welch, W.~J., Mitchell, T.~J., and Wynn, H.~P. (1989).
\newblock Design and analysis of computer experiments.
\newblock {\em Statistical Science\/}, {\bf 4}(4), 409--423.

\bibitem[Santner {\em et~al.}(2003)Santner, Williams, and Notz]{Santner03}
Santner, T.~J., Williams, B.~J., and Notz, W.~I. (2003).
\newblock {\em The Design and Analysis of Computer Experiments\/}.
\newblock Springer Science+Businuess Media New York.

\bibitem[Smith(2013)Smith]{SmithUQ2013}
Smith, R.~C. (2013).
\newblock {\em Uncertainty Quantification: Theory, Implementation, and
  Applications\/}.
\newblock SIAM, Philadelphia.

\bibitem[Sullivan(2015)Sullivan]{SullivanUQ2015}
Sullivan, T. (2015).
\newblock {\em Introduction to Uncertainty Quantification\/}.
\newblock Springer, New York.

\bibitem[Traub and Werschulz(1998)Traub and Werschulz]{Traub1998}
Traub, J.~F. and Werschulz, A.~G. (1998).
\newblock {\em Complexity and Information\/}.
\newblock Cambridge University Press, Cambridge.

\bibitem[Tyagi and Cevher(2014)Tyagi and Cevher]{Tyagi2014}
Tyagi, H. and Cevher, V. (2014).
\newblock Learning non-parametric basis independent models from point queries
  via low-rank methods.
\newblock {\em Applied and Computational Harmonic Analysis\/}, {\bf 37}(3),
  389--412.

\bibitem[Wang and Shan(2006)Wang and Shan]{Wang2006}
Wang, G.~G. and Shan, S. (2006).
\newblock Review of metamodeling techniques in support of engineering design
  optimization.
\newblock {\em Journal of Mechanical Design\/}, {\bf 129}(4), 370--380.

\bibitem[Weisberg(2005)Weisberg]{Weisberg2005}
Weisberg, S. (2005).
\newblock {\em Applied Linear Regression\/}.
\newblock John Wiley \& Sons, Inc., New York, 3rd edition.

\bibitem[Yin {\em et~al.}(2008)Yin, Li, and Cook]{Yin08}
Yin, X., Li, B., and Cook, R.~D. (2008).
\newblock Successive direction extraction for estimating the central subspace
  in a multiple-index regression.
\newblock {\em Journal of Multivariate Analysis\/}, {\bf 99}(8), 1733--1757.

\bibitem[Zhang {\em et~al.}(2017)Zhang, Li, Lin, Zeng, and Wu]{Zhang2017}
Zhang, J., Li, W., Lin, G., Zeng, L., and Wu, L. (2017).
\newblock Efficient evaluation of small failure probability in high-dimensional
  groundwater contaminant transport modeling via a two-stage {Monte Carlo}
  method.
\newblock {\em Water Resources Research\/}, {\bf 53}(3), 1948--1962.

\end{thebibliography}

%
%

\end{document}